\documentclass[10pt]{article}
\usepackage{bbm}
\usepackage{amsmath}
\usepackage{amsthm}
\usepackage{bbm}
\usepackage{amssymb} 
\usepackage{enumitem}
\usepackage{mathrsfs}
\usepackage[dvipsnames]{xcolor}
\usepackage{graphicx}
\usepackage{authblk}
\usepackage[backref=none,hypertexnames=false,colorlinks=true,urlcolor=blue,linkcolor=blue,citecolor=blue]{hyperref}
\usepackage[numbers,comma,square,sort&compress]{natbib}
\usepackage[letterpaper,text={6.5in,9in},centering]{geometry}

\graphicspath{{Figs/}}

\makeatletter\@addtoreset{equation}{section}\makeatother
\renewcommand{\theequation}{\arabic{section}.\arabic{equation}}

\newtheorem{theorem}{Theorem}

\newenvironment{Acknowledgment}%
 {\begin{trivlist}\item[]\textbf{Acknowledgments.}}{\end{trivlist}}

 \newenvironment{Data}%
 {\begin{trivlist}\item[]\textbf{Data Availability Statement.}}{\end{trivlist}}

\begin{document}

\title{
Predicting the emergence of localised dihedral patterns in models for dryland vegetation
}
\author[1]{Dan J. Hill}

\affil[1]{\small Fachrichtung Mathematik,  Universit\"at des Saarlandes, Postfach 151150, 66041 Saarbr\"ucken, Germany}

\date{}
\maketitle

\begin{abstract}\noindent
Localised patterns are often observed in models for dryland vegetation, both as peaks of vegetation in a desert state and as gaps within a vegetated state, known as `fairy circles'. Recent results from radial spatial dynamics show that approximations of localised patterns with dihedral symmetry emerge from a Turing instability in general reaction--diffusion systems, which we apply to several vegetation models. We present a systematic guide for finding such patterns in a given reaction--diffusion model, during which we obtain four key quantities that allow us to predict the qualitative properties of our solutions with minimal analysis.  {  We consider four well-established vegetation models and compute their key predictive quantities, observing that models which possess similar values exhibit qualitatively similar localised patterns; we then complement our results with numerical simulations of various localised states in each model. Here, localised vegetation patches emerge generically from Turing instabilities and act as transient states between uniform and patterned environments, displaying complex dynamics as they evolve over time.}
\end{abstract}

\section{Introduction}\label{s:Intro}

In semi-arid environments, the formation of distinct vegetation patterns has been well documented since their first discovery in sub-Saharan Africa in the 1950s \cite{MacFadyen1950a,MacFadyen1950b}. Sloped terrains exhibit stripes (or bands) of vegetation \cite{Siero2015,Samuelson2019,Deblauwe2011Modulation,Deblauwe2012Migration}, while flat terrains are home to periodic arrangements of spots and gaps, as well as labyrinthine patterns \cite{Gowda2014,Gowda2016,Borgogno2009,Ludwig2005vegpatches,Rietkerk2002}; see \cite{MERON2016Pattern} for a detailed review of the mechanisms that cause vegetation patterns. The environmental role of these structures is still debated: are they a sign of desertification, providing early indications of environmental decline \cite{Rietkerk2004catastrophe}, or are they nature's safety net, allowing an ecosystem to avoid critical tipping points and stave off extinction \cite{Rietkerk2021Tipping}? Such questions have generated significant interest in vegetation patterns over several decades, and this looks set to only increase in the face of our changing climate.

One of the key questions when first studying vegetation patterns is what kind of model to consider---the dynamics of an ecosystem are highly complex, and any attempts to model a particular environment will need to reflect the specific qualities of that environment. Vegetation with vertical roots and strong water uptake are subject to a pattern-forming feedback driven by diffusion, while vegetation with widely spread roots are also subject to nonlocal effects, due to water uptake by the roots \cite{MERON2016Pattern}. Any attempts at modelling these types of environments must also reflect these differences; see \cite{Lejeune2002,Klausmeier1999,vonHardenberg2001,vanderstelt2013} for diffusive PDE models, and \cite{Escaff2015,Meron2007,Tlidi2008,Gilad2004Engineers} for nonlocal models. 

Even with these environment-specific modelling terms there are certain phenomena that appear across many vegetation models, which we call \textit{universal phenomena}; we can better understand these universal phenomena by identifying the mechanisms that cause them. For example, periodic and labyrinthine patterns are known to emerge from bifurcations with non-zero wavenumber, which we refer to as Turing bifurcations, and so most flat terrain vegetation models will undergo such bifurcations for a critical choice of parameters. The presence of a universal phenomenon (such as periodic patterns) implies the presence of a universal mechanism (such as a Turing bifurcation) within a phenomenological model.

A particular phenomenon that we are interested in is the emergence of fully localised patterns, where a patterned state is completely surrounded by a uniform state. A notable example of these are so called `fairy circles', consisting of roughly circular barren patches surrounded by vegetation \cite{Getzin2021Fairy}. Fairy circles have been found to exhibit a consistent wavelength between neighbours, and can be found in fully localised patches or as constituents of a larger periodic structure \cite{Getzin2015Fairy}. Originally observed in the Namib desert \cite{VanRooyen2004Fairy}, fairy circles have since been discovered in Western Australia \cite{Getzin2016FairyCircles}, suggesting that they may be more universal than originally supposed.

\begin{figure}[t!]
    \centering
    \includegraphics[width=0.8\linewidth]{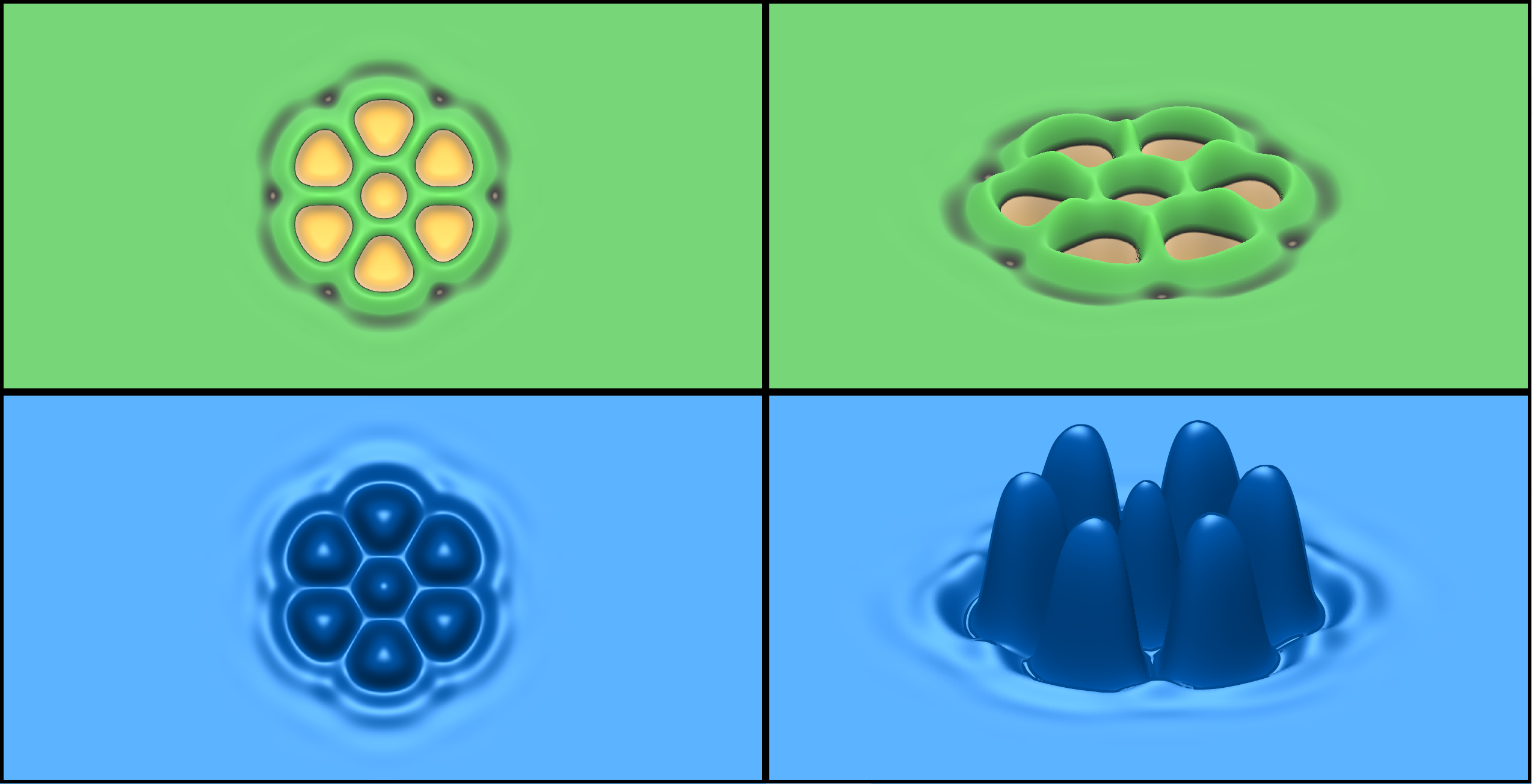}
    \caption{Example of a localised vegetation pattern with hexagonal symmetry, with vegetation density plotted in the top row and water density plotted in the bottom row. The left-hand column provides a top-down view of the pattern, while the right-hand column provides a three-dimensional view. Here the vegetation has a negative polarity and so consists of gaps, whereas the water has a positive polarity and is thus made up of peaks. These images were produced using VisualPDE \cite{Walker2023VisualPDE}. 
    }
    \label{fig:intro-fig}
\end{figure}

While localised patterns are well understood mathematically in one spatial dimension \cite{Woods1999,chapman2009exponential}, this is not the case in two dimensions. There has been recent progress on axisymmetric patterns both in vegetation models \cite{Hill2022Veg,byrnes2023large} and other contexts \cite{lloyd2009localized,mccalla2013spots,mcquighan2014oscillons,Hill2020Localised}, mostly using radial spatial dynamics theory pioneered by Scheel \cite{scheel2003radially}. However, there is very little theory regarding fully localised non-axisymmetric patterns outside of variational methods \cite{Buffoni2022} {  and exponential asymptotics \cite{Kozyreff2013Hexagonal}}. We focus on fully localised patterns with dihedral symmetry---such that they are invariant under rotations of $2\pi/m$, for some $m\in\mathbb{N}$, and one line of reflection---which includes localised hexagons and squares observed in prototypical pattern-forming systems \cite{lloyd2008localized,sakaguchi1997stable}.

We consider a simple class of two-component reaction--diffusion models that cover a wide range of vegetation models, where pattern formation is driven by diffusion and domain-covering patterns emerge from Turing bifurcations. For this general class of models, we show that fully localised patterns with dihedral symmetry emerge universally from a Turing bifurcation. To do this, we first present a step-by-step procedure to determine the local expression of a reaction--diffusion system near a Turing bifurcation, which then allows us to leverage recent rigorous results in \cite{Hill2023DihedralPatch,Hill2024dihedral} regarding the existence of approximate fully localised dihedral patterns. We present theorems for the existence of two classes of localised dihedral patterns: spot A-type patterns, which can appear as either peaks or gaps, and ring-type patterns that emerge from a pitchfork bifurcation. If a pattern mostly consists of peaks we say it has a positive \emph{polarity}, and a negative polarity if it mostly consists of gaps; see Figure~\ref{fig:intro-fig}. 

\begin{figure}[t!]
    \centering
    \includegraphics[width=\linewidth]{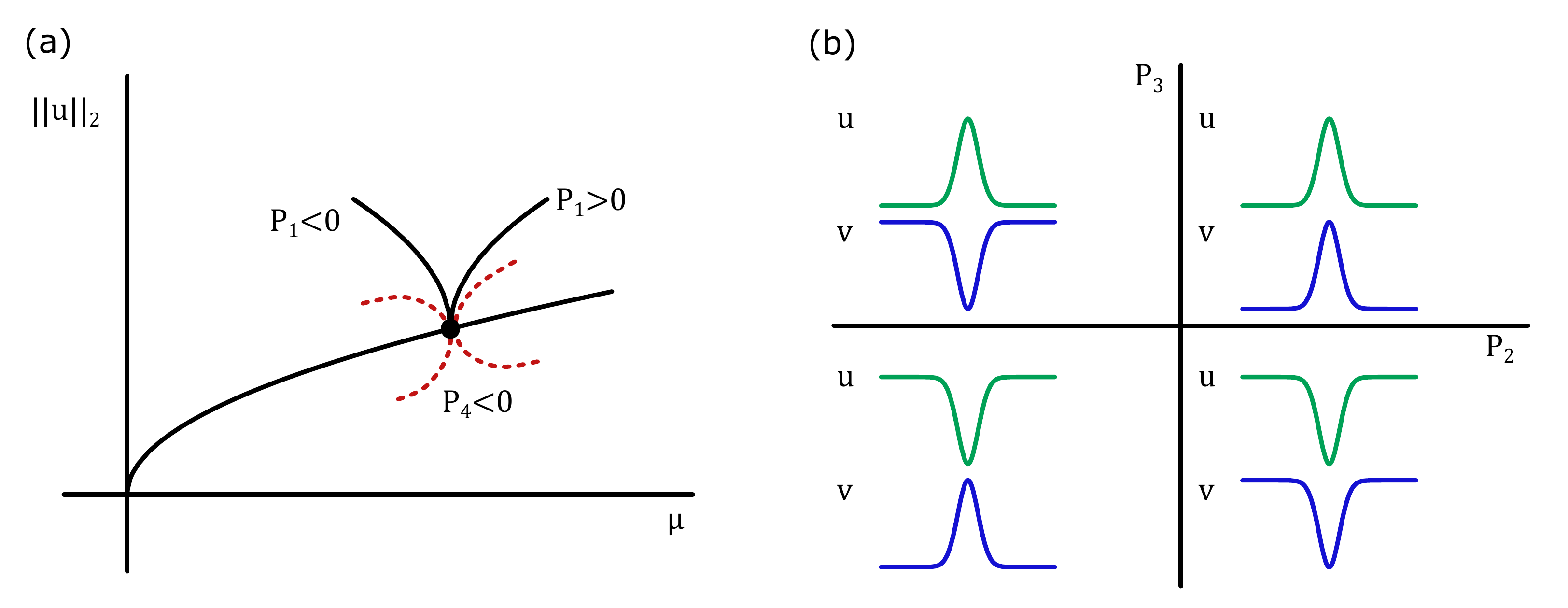}
    \caption{(a) Possible bifurcation curves for localised dihedral patterns. The sign of $P_1$ determines the direction of bifurcation for spot A-type patterns (black, solid) and, if $P_4<0$, ring-type patterns (red, dashed). (b) The different profiles of localised patterns; $P_2$ determines the phase of the solutions, while $P_3$ determines the polarity. Here, $u$ denotes vegetation and $v$ denotes water density.}
    \label{fig:Predictors}
\end{figure}
This process also reveals four key quantities $P_1,\dots,P_4$, which we call \emph{qualitative predictors}, that establish the qualitative behaviour of these localised patterns. In particular,
\begin{enumerate}[label=(\roman*)]
    \item $P_1$ determines the direction of bifurcation from the Turing point; see Figure~\ref{fig:Predictors}(a). If $P_1>0$ then localised solutions bifurcate as the precipitation increases and if $P_1<0$ then localised solutions bifurcate as the precipitation decreases.
    \item $P_2$ determines the relationship between vegetation and water densities {  of spot A-type patterns}; see Figure~\ref{fig:Predictors}(b). If $P_2>0$ then the vegetation and water share the same polarity and if $P_2<0$ then they have opposing polarity. Such vegetation patterns are respectively called \emph{in-phase} or \emph{anti-phase}, both of which have been observed in various semi-arid climates; see \cite{Getzin2016FairyCircles}, for example.
    \item $P_3$ determines the polarity of spot A-type patterns; see Figure~\ref{fig:Predictors}(b). If $P_3>0$ then spot A-type patterns are made up of peaks and if $P_3<0$ then spot A-type patterns are made up of gaps.
    \item $P_4$ determines whether or not ring-type patterns exist. If $P_4<0$ then ring-type patterns bifurcate from the Turing point, as seen in Figure~\ref{fig:Predictors}(a), and if $P_4>0$ then no ring-type patterns emerge.
\end{enumerate}
This work represents the first analytic evidence for fully localised non-axisymmetric patterns in vegetation models, to the author's knowledge, and establishes a guide for quickly analysing such patterns in other two-component reaction--diffusion systems. {  We emphasise that our theoretical results are regarding the \emph{existence} of localised steady-state patterns for reaction--diffusion systems, and so they do not determine the \emph{stability} of such solutions. As such, these patterns instead represent transient states as a uniform vegetated state transitions to a patterned state; we later observe the nonlinear dynamics of this transition in Section~\ref{s:numeric}. Since the time scales associated with vegetation models are often larger than other pattern-forming systems, these transient states are more likely to be observed in the field.}

We find that spot A-type patterns are highly universal, emerging generically from a Turing bifurcation, whereas ring-type patterns only emerge in narrow parameter regions. In each of the vegetation models we consider spot A-type patterns grow in width over time, invading the surrounding uniform state. As such, these localised patterns can represent a degradation of a uniform vegetated state resulting in a significant change to the environment, such as desertification. 

The rest of the paper is structured as follows: in Section~\ref{s:Theory-deriv} we derive the local expression for a general two-component reaction--diffusion system near a Turing instability, presenting a step-by-step process that can be easily automated. In Section~\ref{s:Theory-veg} we then introduce four notable reaction--diffusion models for dryland vegetation, and compute their qualitative predictors in each case. Finally, in Section~\ref{s:numeric} we use exponential time-steppers to numerically simulate three examples of localised dihedral patterns in each vegetation model, comparing the similarities and differences between different models.

\section{Theory}\label{s:Theory}
We begin this section by taking a general two-component reaction--diffusion model and deriving its local representation in the neighbourhood of a Turing instability. This allows us to utilise abstract theorems to predict the existence of solutions in various vegetation models. Throughout this procedure we identify several notable quantities which describe the qualitative behaviour of our localised solutions; computing these key quantities then allows us to make heuristic predictions regarding the emergence of localised patterns in each specific reaction--diffusion model. Each stage of the procedure is reduced to solving an algebraic condition, thus the entire process can be easily automated for rapid predictions regarding the emergence of localised patterns in two-component reaction--diffusion models.

Following this, we introduce four vegetation models that belong to our general class of two-component reaction--diffusion equations: the Klausmeier--Gray--Scott, logistic Klausmeier, NFC--Gilad, and von Hardenberg models. For each model we compute our key quantities and predict the behaviour of localised solutions bifurcating from each Turing point, we will then compare these predictions with the numerical results presented in Section~\ref{s:numeric}.
\subsection{Derivation}\label{s:Theory-deriv}
In general, a reaction--diffusion equation has the form
\begin{equation}\label{eqn:RD-N}
    \partial_{t}\mathbf{u} = \mathbf{D}\Delta\mathbf{u} + \mathbf{F}(\mathbf{u};\mu).
\end{equation}
Here, $\mathbf{u}(t,\mathbf{x})\in\mathbb{R}^{N}$ denotes $N$ coupled variables, with temporal and spatial coordinates $t\in\mathbb{R}$ and $\mathbf{x}\in\mathbb{R}^{2}$, respectively, and $\mu\in\mathbb{R}^{p}$ denotes $p$ parameters of the system, for some $N,p\in\mathbb{N}$. The two-dimensional Laplacian $\Delta:=\partial_x^2 + \partial_y^2$ encodes the spatial diffusion of $\mathbf{u}$, and the $N\times N$ dimensional square matrix $\mathbf{D}\in\mathbb{R}^{N\times N}$ defines the diffusion coupling of the system (i.e., how the spatial diffusion of each variable affects the overall temporal growth of $\mathbf{u}$). Finally, the function $\mathbf{F}\in\mathbb{R}^{N}$ accounts for all local reactions in the model.

For this present work we assume that $\mathbf{u}(t,\mathbf{x})$ is made up of two variables $u(t,\mathbf{x}),v(t,\mathbf{x})\geq 0$ which represent the densities of the vegetation biomass and soil water, respectively, and that the diffusion matrix $\mathbf{D}\in\mathbb{R}^{2\times2}$ is lower triangular, such that the rate of growth of the biomass is independent of the diffusive rate of the water. This assumption models the effect of vegetation on water diffusion through root suction. We also assume that $p-1$ parameters are fixed, such that any bifurcations are driven by a single parameter $\mu\in\mathbb{R}$, which we call the bifurcation parameter (in vegetation models $\mu$ is usually the precipitation rate). 

As such, the remainder of this work will focus on general reaction--diffusion systems of the form
\begin{equation}\label{eqn:RD-gen}
    \begin{split}
        u_t &= \Delta u - \hat{f}(u,v;\mu),\\
        v_t &= D_v\Delta (v-\beta u) - \hat{g}(u,v;\mu).\\
    \end{split}
\end{equation}
We have defined the system \eqref{eqn:RD-gen} so that subsequent analysis will be slightly easier; a more general system like
\begin{equation*}
    \begin{split}
        u_t &= D_{uu}\Delta u  + \tilde{f}(u,v;\mu),\\
        v_t &= D_{vv}\Delta v + D_{uv}\Delta u + \tilde{g}(u,v;\mu)\\
    \end{split}
\end{equation*}
could equivalently be written in the form of \eqref{eqn:RD-gen} by rescaling the coordinate system $(t,\mathbf{x})$ with respect to $D_{uu}$ and then defining parameters $D_v,\beta$ and functions $\hat{f},\hat{g}$ accordingly. We are interested in steady states of the system \eqref{eqn:RD-gen} and so, for the rest of this section, we will consider the following equations
\begin{equation}\label{eqn:RD-Stationary}
    \begin{split}
        0 &= \Delta u - f(u,v;\mu),\\
        0 &= \Delta v - g(u,v;\mu),\\
    \end{split}
\end{equation}
where we have defined
\begin{equation*}
    \begin{split}
        f(u,v;\mu) &:= \hat{f}(u,v;\mu),\\
        g(u,v;\mu) &:= \frac{1}{D_v}\hat{g}(u,v;\mu) + \beta\hat{f}(u,v;\mu).\\
    \end{split}
\end{equation*}
In Section~\ref{s:numeric}, we will then return to the full time-dependent PDE system \eqref{eqn:RD-gen} for our numerical investigations. 

The aim for this section is to use abstract theory to predict the emergence of fully localised dihedral patterns in \eqref{eqn:RD-Stationary} near a Turing instability. By Turing instability, we mean a parameter value at which a uniform steady state loses stability with respect to perturbations with nonzero wave number $k>0$; sometimes the term Turing instability is reserved for instabilities caused by changes in the diffusion coefficient $\mathbf{D}$, but we will not use that definition in this paper. To do this, we present the following procedure:
\begin{enumerate}[label=(\S2.1.\arabic*)]
    \item Identify uniform steady states of the system \eqref{eqn:RD-Stationary}.
    \item Compute the spatial eigenvalues of the system \eqref{eqn:RD-Stationary} at the uniform steady state and identify a critical parameter value $\mu_*(k)$ such that \eqref{eqn:RD-Stationary} undergoes a Turing bifurcation with wave number $k>0$.
    \item Express the system \eqref{eqn:RD-Stationary} as a Taylor expansion about the uniform steady state at the Turing point $\mu=\mu_*(k)$.
    \item Use abstract existence results from \cite{Hill2023DihedralPatch,Hill2024dihedral} for localised approximate dihedral patterns near a Turing instability.
\end{enumerate}
Once we have transformed the system \eqref{eqn:RD-Stationary} into local coordinates in the neighbourhood of a Turing instability, we identify four key quantities $P_1,\dots,P_4$ that help us predict various properties of localised dihedral patterns bifurcating from the Turing point $(\S2.1.5)$.

\subsubsection{Finding uniform steady states}\label{ss:1}
The first task is to identify uniform solutions of the steady state equation \eqref{eqn:RD-Stationary}. That is, we look for solutions that are independent of the spatial coordinate $\mathbf{x}$, such that $\Delta u = \Delta v = 0$. Then, the system \eqref{eqn:RD-Stationary} becomes
\begin{equation}\label{eqn:uniform}
    \begin{split}
        0 &= f(u,v;\mu),\qquad\qquad
        0 = g(u,v;\mu)\\
    \end{split}
\end{equation}
for which we want to find solutions $u = u_*(\mu)$, $v=v_*(\mu)$. In general, while these equations can be difficult to solve by hand, especially with additional fixed parameters in a given model, they are very tractable for numerical solvers. For the vegetation models presented later in this section, the problem reduces to solving a cubic-order polynomial, and so one would typically expect to find between one and three uniform solutions to \eqref{eqn:uniform} for any fixed $\mu\in\mathbb{R}$. 

Having identified some uniform steady state $(u,v;\mu) = \mathcal{U}_*(\mu):=(u_*(\mu), v_*(\mu); \mu)$, we can restrict our analysis to a small neighbourhood of $\mathcal{U}_*(\mu)$ by taking the Taylor series of $f$ and $g$. {  We note that, in order to perform weakly nonlinear analysis in Section 2.1.4, we require explicit formulae for the quadratic and cubic nonlinear terms in the Taylor expansion of $f$ and $g$.} Thus, we define $u(\mathbf{x}) = u_*(\mu) + U(\mathbf{x})$ and $v(\mathbf{x}) = v_*(\mu) + V(\mathbf{x})$, with $|U|,|V|\ll1$, so that
\begin{equation}\label{eqn:Taylor-mu}
    \begin{split}
        0 &= \Delta \mathbf{U} - \mathbf{M}(\mu)\mathbf{U} - \mathbf{F}_2(\mathbf{U};\mu)  - \mathbf{F}_3(\mathbf{U};\mu) + \mathcal{O}(|\mathbf{U}|^4),\\
    \end{split}
\end{equation}
where $\mathbf{U}:=(U,V)^T$ and 
\begin{equation}
    \begin{split}
        \mathbf{M}(\mu) &:= \begin{pmatrix} \partial_{u}f(\mathcal{U}_{*}(\mu)) & \partial_{v}f(\mathcal{U}_{*}(\mu)) \\ \partial_{u}g(\mathcal{U}_{*}(\mu)) & \partial_{v}g(\mathcal{U}_{*}(\mu)) \end{pmatrix},\\
        \mathbf{F}_{2}(\mathbf{U};\mu) &:= \frac{1}{2}\begin{pmatrix} \partial_{u}^{2}f(\mathcal{U}_{*}(\mu))U^2 +  2\partial_{u}\partial_{v}f(\mathcal{U}_{*}(\mu))U V + \partial_{v}^{2}f(\mathcal{U}_{*}(\mu))V^2 \\ \partial_{u}^{2}g(\mathcal{U}_{*}(\mu))U^2 + 2\partial_{u}\partial_{v}g(\mathcal{U}_{*}(\mu))U V + \partial_{v}^{2}g(\mathcal{U}_{*}(\mu))V^2 \end{pmatrix},\\
        \mathbf{F}_{3}(\mathbf{U};\mu) &:= \frac{1}{6}\begin{pmatrix} \partial_{u}^{3}f(\mathcal{U}_{*}(\mu))U^3 + 3\partial_{u}^{2}\partial_{v}f(\mathcal{U}_{*}(\mu))U^2 V +  3\partial_{u}\partial_{v}^{2}f(\mathcal{U}_{*}(\mu)) U V^2 + \partial_{v}^{3}f(\mathcal{U}_{*}(\mu))V^3 \\ \partial_{u}^{3}g(\mathcal{U}_{*}(\mu))U^3 + 3\partial_{u}^{2}\partial_{v} g(\mathcal{U}_{*}(\mu))U^2 V +  3\partial_{u}\partial_{v}^{2}g(\mathcal{U}_{*}(\mu))U V^2 + \partial_{v}^{3}g(\mathcal{U}_{*}(\mu))V^3 \end{pmatrix}.\\
    \end{split}
\end{equation}
Here, the matrix $\mathbf{M}(\mu)$ is the linear operator of \eqref{eqn:Taylor-mu}, $\mathbf{F}_{2}$ contains all quadratic (in $\mathbf{U}$) nonlinearities of \eqref{eqn:Taylor-mu}, and $\mathbf{F}_{3}$ contains all cubic (in $\mathbf{U}$) nonlinearities of \eqref{eqn:Taylor-mu}. Note that we do not truncate the Taylor series at cubic order---the $\mathcal{O}(|\mathbf{U}|^4)$ terms remain in our problem but do not affect our existence results, and so we choose not to write them down explicitly. 

\subsubsection{Spatial eigenvalue analysis}\label{ss:2}
We now wish to compute the spatial eigenvalues of the linear operator for the system \eqref{eqn:Taylor-mu}; in particular, for fixed $\mu\in\mathbb{R}$ we wish to find $\lambda \in\mathbb{C}$ such that 
\begin{equation*}
    \begin{split}
        0 &= \lambda \mathbf{U} - \mathbf{M}(\mu)\mathbf{U}.\\
    \end{split}
\end{equation*}
This is equivalent to solving the eigenvalue problem of $\mathbf{M}(\mu)$ for $\lambda$ and so we need to solve $\sigma(\lambda;\mu)=0$, where
\begin{equation}
\sigma(\lambda;\mu) := \det\,\left(\lambda\mathbbm{1} - \mathbf{M}(\mu)\right).
\end{equation}
{  We briefly comment on how $\sigma(\lambda;\mu)$ connects to the  standard dispersion relation found in one-dimensional Turing analysis. Taking the original time-dependent system \eqref{eqn:RD-gen} for $\mathbf{u}\in\mathbb{R}^{2}$, restricted to one spatial dimension, and introducing the ansatz $\mathbf{u}(t,x) = \mathbf{u}_* + \hat{\mathbf{u}}\,\mathrm{e}^{\omega t}\,\mathrm{e}^{\mathrm{i}kx} + c.c.$ with $\mathbf{u}_*$ a uniform steady state, we obtain the following dispersion relation
\begin{equation*}
    d(\omega,k;\mu) = \mathrm{det}\,\left(-k^2\mathbf{D} + \partial_{\mathbf{u}}\mathbf{F}(\mathbf{u}_*;\mu) - \omega\mathbbm{1}\right) = 0.
\end{equation*}
Then the condition $\sigma(\lambda;\mu)=0$ corresponds to the critical value $\omega=0$,
 such that the uniform state $\mathbf{u}_*$ changes stability with respect to perturbations with spatial eigenvalues $\pm\sqrt{\lambda}$.} Using the above definition of $\mathbf{M}(\mu)$, we find 
\begin{equation*}\begin{split}
    \sigma(\lambda;\mu) &= \left(\lambda -\partial_{u}f(\mathcal{U}_{*}(\mu)) \right)\left(\lambda -\partial_{v}g(\mathcal{U}_{*}(\mu)) \right) - \partial_{v}f(\mathcal{U}_{*}(\mu))\;\partial_{u}g(\mathcal{U}_{*}(\mu)),\\
    &= \left(\lambda - \tfrac{\left(\partial_{u}f(\mathcal{U}_{*}(\mu)) + \partial_{v}g(\mathcal{U}_{*}(\mu)) \right)}{2} \right)^2 -  \tfrac{4\partial_{v}f(\mathcal{U}_{*}(\mu))\;\partial_{u}g(\mathcal{U}_{*}(\mu)) +\left(\partial_{u}f(\mathcal{U}_{*}(\mu)) - \partial_{v}g(\mathcal{U}_{*}(\mu)) \right)^2}{4}.
\end{split}\end{equation*}

\begin{figure}[t!]
    \centering
    \includegraphics[width=\linewidth]{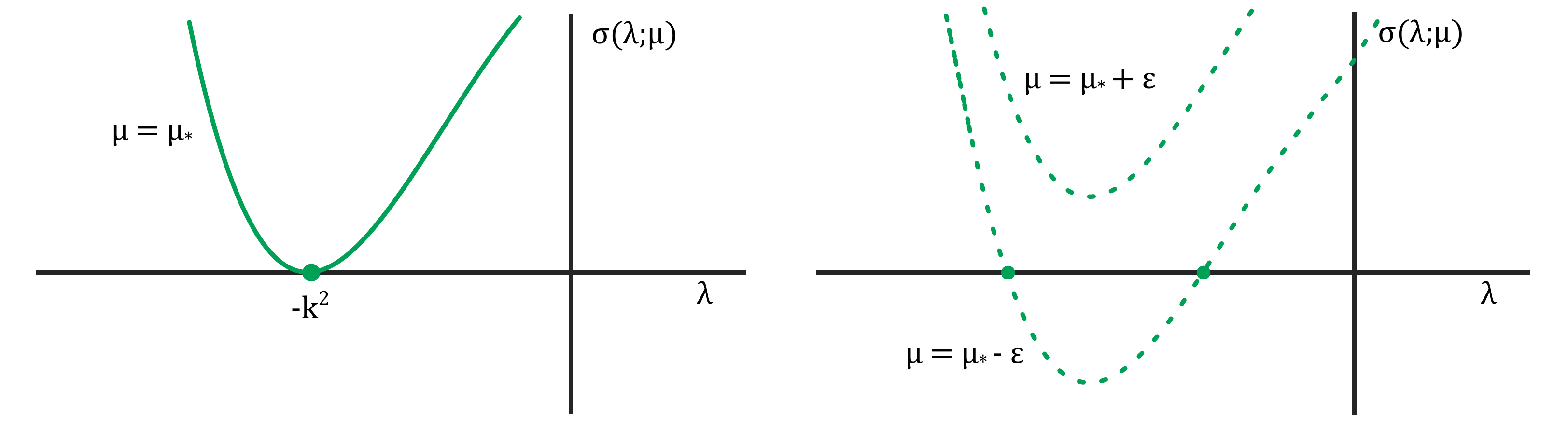}
    \caption{Schematic diagram for the roots of $\sigma(\lambda;\mu_*)$, including the change in $\sigma(\lambda;\mu)$ as $\mu$ varies.}
    \label{fig:disp}
\end{figure}

The point $\mu=\mu_*\in\mathbb{R}$ defines a Turing instability if $\sigma(-k^2;\mu_*) = 0$ for some $k\in\mathbb{R}$, $\sigma(\lambda;\mu_*+\varepsilon)$ has no real roots, and $\sigma(\lambda;\mu_*-\varepsilon)$ has two distinct negative roots for $0<|\varepsilon|\ll1$; see Figure~\ref{fig:disp}. {  Note that we have adopted the convention that $\sigma(\lambda;\mu_*+\varepsilon)$ has no real roots but we have not assumed the sign of $\varepsilon$, and so $\varepsilon$ could be positive or negative.} Thus, $\mu=\mu_*$ is a Turing point if
\begin{equation*}\begin{split}
    \sigma(\lambda;\mu_*) &= \left(\lambda - \tfrac{\left(\partial_{u}f(\mathcal{U}_{*}(\mu_*)) + \partial_{v}g(\mathcal{U}_{*}(\mu_*)) \right)}{2} \right)^2 -  \tfrac{4\partial_{v}f(\mathcal{U}_{*}(\mu_*))\;\partial_{u}g(\mathcal{U}_{*}(\mu_*)) +\left(\partial_{u}f(\mathcal{U}_{*}(\mu_*)) - \partial_{v}g(\mathcal{U}_{*}(\mu_*)) \right)^2}{4}
\end{split}\end{equation*}
possesses repeated real roots, or equivalently
\begin{equation}\begin{split}
    0 &= 4\partial_{v}f(\mathcal{U}_{*}(\mu_*))\;\partial_{u}g(\mathcal{U}_{*}(\mu_*)) +\left(\partial_{u}f(\mathcal{U}_{*}(\mu_*)) - \partial_{v}g(\mathcal{U}_{*}(\mu_*)) \right)^2.
\end{split}\end{equation}
Supposing that $\mu_*$ satisfies the above condition, the wave number $k\in\mathbb{R}$ is then determined to be
\begin{equation}\label{defn:k}\begin{split}
    k^2 &= -\left(\tfrac{\partial_{u}f(\mathcal{U}_{*}(\mu_*)) + \partial_{v}g(\mathcal{U}_{*}(\mu_*)) }{2}\right).
\end{split}\end{equation}
We note that $k^2<0$ would imply that $\mu_*$ is not a Turing point, and so we proceed under the assumption that the right-hand side of \eqref{defn:k} is positive.

\subsubsection{Local coordinates near a Turing instability}\label{ss:4}
We centre our system \eqref{eqn:Taylor-mu} at the point $\mu=\mu_*$ where the uniform state $\mathcal{U}_*(\mu)$ undergoes a Turing instability. To do this, we define $\mu = \mu_* + \varepsilon$, for $0<|\varepsilon|\ll1$, such that \eqref{eqn:Taylor-mu} becomes
\begin{equation}\label{eqn:Taylor-eps}
    \begin{split}
        0 &= \Delta \mathbf{U} - \mathbf{M}_{1}\mathbf{U} -\varepsilon\mathbf{M}_{2}\mathbf{U} - \mathbf{Q}(\mathbf{U},\mathbf{U}) - \mathbf{C}(\mathbf{U},\mathbf{U},\mathbf{U}) + \mathcal{O}\left(|\mathbf{U}|^4 + |\varepsilon|^2\,|\mathbf{U}| + |\varepsilon|\,|\mathbf{U}|^2\right).\\
    \end{split}
\end{equation}
Here, we have introduced new objects $\mathbf{M}_{1},\mathbf{M}_{2},\mathbf{Q},\mathbf{C}$ to bring our system in line with those considered in \cite{Hill2023DihedralPatch,Hill2024dihedral}. In particular, we define $\mathbf{M}_{1} := \mathbf{M}(\mu_{*})$ and $\mathbf{M}_{2} := \mathbf{M}'(\mu_{*})$ as the first terms in the Taylor expansion of $\mathbf{M}(\mu_* + \varepsilon)$. For the nonlinear terms $\mathbf{F}_{2}$, $\mathbf{F}_{3}$ we have defined
\begin{equation}
    \begin{split}
        \mathbf{Q}(\mathbf{X},\mathbf{Y}) &:= \tfrac{1}{2}\left[\mathbf{F}_{2}(\mathbf{X}+\mathbf{Y};\mu_*) - \mathbf{F}_{2}(\mathbf{X};\mu_*) - \mathbf{F}_{2}(\mathbf{Y};\mu_*)\right], \\
        \mathbf{C}(\mathbf{X},\mathbf{Y},\mathbf{Z}) &:=\tfrac{1}{6}\left[\mathbf{F}_{3}(\mathbf{X}+\mathbf{Y}+\mathbf{Z};\mu_*) - \mathbf{F}_{3}(\mathbf{X}+\mathbf{Y};\mu_*) - \mathbf{F}_{3}(\mathbf{X}+\mathbf{Z};\mu_*) - \mathbf{F}_{3}(\mathbf{Y}+\mathbf{Z};\mu_*)\right.\\
        &\qquad\qquad \left. + \mathbf{F}_{3}(\mathbf{X};\mu_*) + \mathbf{F}_{3}(\mathbf{Y};\mu_*) + \mathbf{F}_{3}(\mathbf{Z};\mu_*)\right]\\
    \end{split}
\end{equation}
so that $\mathbf{Q}(\mathbf{U},\mathbf{U}) = \mathbf{F}_{2}(\mathbf{U};\mu_*)$, $\mathbf{C}(\mathbf{U},\mathbf{U},\mathbf{U}) = \mathbf{F}_{3}(\mathbf{U};\mu_*)$, and $\mathbf{Q}$ and $\mathbf{C}$ are respective symmetric bilinear and trilinear functions, i.e. they satisfy
\begin{equation*}\begin{aligned}
    \mathbf{Q}(\alpha\mathbf{X}_1 + \beta\mathbf{X}_2,\mathbf{Y}) ={}& \alpha\mathbf{Q}(\mathbf{X}_1,\mathbf{Y}) + \beta\mathbf{Q}(\mathbf{X}_2,\mathbf{Y}), &\quad \mathbf{Q}(\mathbf{X},\mathbf{Y}) ={}& \mathbf{Q}(\mathbf{Y},\mathbf{X}),\\
    \mathbf{C}(\alpha\mathbf{X}_1 + \beta\mathbf{X}_2,\mathbf{Y},\mathbf{Z}) ={}& \alpha\mathbf{C}(\mathbf{X}_1,\mathbf{Y},\mathbf{Z}) + \beta\mathbf{C}(\mathbf{X}_2,\mathbf{Y},\mathbf{Z}), &\quad \mathbf{C}(\mathbf{X},\mathbf{Y},\mathbf{Z}) ={}& \mathbf{C}(\mathbf{Y},\mathbf{X},\mathbf{Z}) = \mathbf{C}(\mathbf{Z},\mathbf{Y},\mathbf{X}),\\
\end{aligned}
\end{equation*}
for any $\alpha,\beta\in\mathbb{R}$ and $\mathbf{X}_{1},\mathbf{X}_{2},\mathbf{X},\mathbf{Y},\mathbf{Z}\in\mathbb{R}^{2}$. 

There are a number of quantities in \eqref{eqn:Taylor-eps} that play a part in the existence theorems of \cite{Hill2023DihedralPatch,Hill2024dihedral}; these are mostly not unique to the existence of localised dihedral patterns and so we present them here before starting the next part of our procedure. We begin by defining the following vectors
\begin{equation}\begin{aligned}
\hat{U}_0 &= \begin{pmatrix}
    \partial_{v}f(\mathcal{U}_{*}(\mu_*)) \\ -\left(k^2 + \partial_{u}f(\mathcal{U}_{*}(\mu_*))\right)
\end{pmatrix},& \qquad \hat{U}_1 &= \begin{pmatrix}
     0 \\ k^2
\end{pmatrix}, \\
 \hat{U}^*_0 &= \frac{1}{\partial_{v}f(\mathcal{U}_{*}(\mu_*))}\begin{pmatrix}
     1 \\ 0
\end{pmatrix}, & \qquad \hat{U}^*_1 &= \frac{1}{k^2 \partial_{v}f(\mathcal{U}_{*}(\mu_*))}\begin{pmatrix}
       \left(k^2 + \partial_{u}f(\mathcal{U}_{*}(\mu_*))\right) \\ \partial_{v}f(\mathcal{U}_{*}(\mu_*))    
\end{pmatrix}
\end{aligned}\end{equation}
so that 
\begin{equation*}
    \mathbf{M}_1 \hat{U}_0 = -k^2 \hat{U}_0,\qquad\qquad \mathbf{M}_1 \hat{U}_1 = -k^2 \hat{U}_1 + k^2 \hat{U}_0, \qquad\qquad \hat{U}_{i}^{*}\cdot\hat{U}_{j} = \delta_{i,j}.
\end{equation*}
Then, $\hat{U}_0$, $\hat{U}_1$ are generalised eigenvectors of $\mathbf{M}_1$ such that setting $\mathbf{U} = A\hat{U}_0 + B\hat{U}_1$ transforms \eqref{eqn:Taylor-eps} into a canonical form for a Turing bifurcation. The vectors $\hat{U}_{0}^{*},\hat{U}_{1}^{*}$ form a dual basis for $\hat{U}_{0},\hat{U}_{1}$, allowing us to project \eqref{eqn:Taylor-eps} onto the respective coordinates $A,B$. 

We now introduce the following constants
\begin{equation}
    \begin{split}
        c_0:={}&\hat{U}_{1}^{*}\cdot\left(-\tfrac{1}{4}\mathbf{M}_2 \hat{U}_0\right),\\
     c_3 :={}& -\left(\tfrac{5}{6}\left(\hat{U}_0^*\cdot { Q_{00}}\right) + \tfrac{5}{6}\left(\hat{U}_1^*\cdot { Q_{01}}\right) + \tfrac{19}{18}\left(\hat{U}_1^*\cdot { Q_{00}}\right)\right) \left(\hat{U}_1^*\cdot { Q_{00}}\right) - \tfrac{3}{4} \left(\hat{U}_1^*\cdot { C_{000}}\right),
    \end{split}
\end{equation}
where we have written ${ Q_{ij}}:=\mathbf{Q}(\hat{U}_{i},\hat{U}_{j})$, ${ C_{ijk}}:=\mathbf{C}(\hat{U}_{i},\hat{U}_{j},\hat{U}_{k})$. The constant $c_0$ determines the direction of bifurcation from the Turing point---we will see this explicitly in the upcoming theorems---while the constant $c_3$ determines whether domain-covering stripes (otherwise known as rolls) emerge via a subcritical or supercritical pitchfork bifurcation \cite{Iooss1993Homoclinic}. This second definition may seem irrelevant in our current study; however, it was previously determined that one dimensional localised patterns only bifurcate from the Turing point if rolls bifurcate subcritically \cite{Woods1999}. Finally, we introduce the constant $\gamma$, where
\begin{equation}
    \gamma := \hat{U}_{1}^{*}\cdot { Q_{00}}.
\end{equation}
In the study of axisymmetric patterns, the sign of $\gamma$ determines whether localised spot solutions have an elevation (a peak) or a depression (a gap) at their centre \cite{lloyd2009localized,Hill2020Localised}. We note that the quantities $c_0,c_3,\gamma$ and the eigenvector $\hat{U}_0$ will be important when defining our qualitative predictors $P_1,\dots,P_4$ later in this section. We are now equipped to study the existence theorems for localised approximate dihedral patterns seen in \cite{Hill2023DihedralPatch,Hill2024dihedral}.

\subsubsection{Approximate localised dihedral patterns}\label{ss:5}
We consider fully localised patterns with $\mathbb{D}_{m}$ dihedral symmetry, so that 
\begin{equation}
    \mathbf{U}(r,\theta) = \mathbf{U}(r,\theta + \tfrac{2\pi}{m}) = \mathbf{U}(r,-\theta),
\end{equation}
where $r,\theta$ are the standard polar coordinates for $\mathbb{R}^{2}$. We employ a Galerkin scheme for approximating $\mathbb{D}_{m}$ solutions of the system \eqref{eqn:Taylor-eps}, which we now detail. Taking a truncated Fourier expansion
\begin{equation}
    \mathbf{U}(r) = \sum_{n=-N}^{N} \mathbf{U}_{|n|}(r) \cos(m n \theta),
\end{equation}
for a fixed choice of $N\geq0$ and projecting onto each $\cos(mn\theta)$, \eqref{eqn:Taylor-eps} becomes
\begin{equation}\label{eqn:Galerk}
    \begin{split}
        0 &= \Delta_{n} \mathbf{U}_{n} - \mathbf{M}_{1}\mathbf{U}_{n} -\varepsilon\mathbf{M}_{2}\mathbf{U}_{n} - \sum_{\substack{i+j=n\\|i|,|j|\leq N}}\mathbf{Q}(\mathbf{U}_{|i|},\mathbf{U}_{|j|}) - \sum_{\substack{i+j+k=n\\|i|,|j|,|k|\leq N}}\mathbf{C}(\mathbf{U}_{|i|},\mathbf{U}_{|j|},\mathbf{U}_{|k|}) + h.o.t.s,\\
    \end{split}
\end{equation}
for each $n\in[0,N]$, where $h.o.t.s$ denotes the higher order terms seen in \eqref{eqn:Taylor-eps} and $\Delta_n := (\partial^2_{r} + \tfrac{1}{r}\partial_r - \tfrac{(mn)^2}{r^2})$ is the Laplacian operator $\Delta$ when applied to the $n^\text{th}$ Fourier mode. Note that our introduction of $\mathbf{Q}$ and $\mathbf{C}$ has provided a much cleaner expression of the nonlinear terms in \eqref{eqn:Galerk} than if we used $\mathbf{F}_{2}$ and $\mathbf{F}_{3}$. 

We now write down existence theorems for localised solutions to \eqref{eqn:Galerk} originally proven in \cite{Hill2023DihedralPatch,Hill2024dihedral}. We shorten our statement of the theorems for the sake of brevity and readability, presenting only the key points relevant to this work. First, we consider spot A-type solutions, whose existence was the focus of \cite{Hill2023DihedralPatch}:

\begin{theorem}[\cite{Hill2023DihedralPatch}]\label{thm:Patch} 
	Fix $m\in \mathbb{N}$, $N\in\mathbb{N}_{0}${ . If $\gamma\neq0$, then there exists some $r_0>0$ such that the Galerkin system \eqref{eqn:Galerk} has an approximate localised dihedral solution $\mathbf{U}_{A}(r,\theta)$ for all $r\in[0,\infty)$, $\theta\in[0,2\pi)$ in the region $0<c_0\varepsilon\ll1$, such that $\mathbf{U}_{A}$ has a leading-order expansion}
\begin{equation}\label{RadialProfile}
\mathbf{U}_{A}(r,\theta) = \frac{2\sqrt{3} k}{\gamma } \;(c_0 \varepsilon)^{\frac{1}{2}}\sum_{n=-N}^{N} a_{|n|}\,J_{|m n|}(k r) \cos(m n \theta)\hat{U}_0 + \mathcal{O}(|\varepsilon|)
\end{equation} 
for $r\in[0,r_0]$ and $|\mathbf{U}_{A}{ (r,\theta)}|\to0$ exponentially fast as $r\to\infty$. Here $J_{n}$ is the $n^\mathrm{th}$ order Bessel function of the first kind and the constants $\{a_{n}\}_{n=0}^{N}$ are nondegenerate solutions of the quadratic matching condition 
		\begin{equation}\label{MatchEq}
    a_{n} = 2\sum_{j=1}^{N-n} \cos\left(\frac{m\pi(n-j)}{3}\right) a_{j} a_{n+j} + \sum_{j=0}^{n} \cos\left(\frac{m\pi(n-2j)}{3}\right) a_{j}a_{n-j},
\end{equation}
    for each $n\in[0,N]$.
\end{theorem}
Setting $N=0$ recovers the axisymmetric spot A solution found in \cite{lloyd2009localized} for the Swift--Hohenberg equation, which is why we refer to these solutions as spot A-type patterns. The proof of Theorem~\ref{thm:Patch} requires tools from radial spatial dynamics, pioneered by Scheel \cite{scheel2003radially}, and can be found in Section 4 of \cite{Hill2023DihedralPatch}. Solving the algebraic matching condition \eqref{MatchEq} then provides excellent initial approximations of fully localised dihedral patterns for numerical study; we will use an initial guess of \eqref{RadialProfile} with solutions from \eqref{MatchEq} as the starting point for our numerical simulations in Section~\ref{s:numeric}. 

Beyond the spot A-type dihedral patterns found in \cite{Hill2023DihedralPatch}, we can also find localised dihedral ring patterns for the system \eqref{eqn:Galerk}; the existence of these patterns was the focus of \cite{Hill2024dihedral}, resulting in the following theorem.
\begin{theorem}[\cite{Hill2024dihedral}]\label{thm:Ring} 
	Fix $m\in \mathbb{N}$, $N\in\mathbb{N}_{0}${ . If $c_3<0$, then there exists some $r_0>0$ such that the Galerkin system \eqref{eqn:Galerk} has an approximate localised dihedral solution $\mathbf{U}_{R}(r,\theta)$ for all $r\in[0,\infty)$, $\theta\in[0,2\pi)$ in the region $0<c_0\varepsilon\ll1$, such that $\mathbf{U}_{R}$ has a leading-order expansion}
\begin{equation}\label{RadialProfile-ring}
\mathbf{U}_{R}(r,\theta) = C_{R}(c_0\varepsilon)^{\frac{3}{4}}\sum_{n=-N}^{N} b_{|n|}\,\left[k r J_{|m n + 1|}(k r)\hat{U}_0 + 2 J_{|m n|}(k r)\hat{U}_1 \right]\cos(m n \theta) + \mathcal{O}(|\varepsilon|)
\end{equation} 
for $r\in[0,r_0]$ and $|\mathbf{U}_{R}({ r,\theta})|\to0$ exponentially fast as $r\to\infty$. Here $C_{R}>0$ is a fixed constant, $J_{n}$ is the $n^\mathrm{th}$ order Bessel function of the first kind, and the constants $\{b_{n}\}_{n=0}^{N}$ are nondegenerate solutions of the cubic matching condition 
		\begin{equation}\label{MatchEq-cubic}
    b_{n} = \sum_{\substack{i+j+k=n\\|i|,|j|,|k|\leq N}} (-1)^{\frac{m(|i|+|j|-|k|-n)}{2}} b_{|i|} b_{|j|} b_{|k|},
\end{equation}
    for each $n\in[0,N]$.
\end{theorem}
Setting $N=0$ recovers the axisymmetric ring solution found in \cite{lloyd2009localized} for the Swift--Hohenberg equation, and so we refer to these solutions as ring-type patterns. The proof of Theorem~\ref{thm:Ring} follows in a similar way to Theorem~\ref{thm:Patch}, but with more delicate analysis regarding the exponential decay of solutions for large values of $r$. As discussed earlier, these ring-type patterns only emerge when bifurcating stripes are unstable; in particular, their far-field profile is approximated by the localised solution of the cubic nonautonomous Ginzburg--Landau equation
\begin{equation}\label{eqn:Ginz}
    \left(\frac{\mathrm{d}}{\mathrm{d}s} + \frac{1}{2s}\right)^{2}q(s) = q(s) + c_3 q(s)^3,
\end{equation}
where $s:=(c_0\varepsilon)^{\frac{1}{2}}r$ is a rescaled radial coordinate in the far-field. An exponentially decaying solution to \eqref{eqn:Ginz} only exists when $c_3<0$, as proven in \cite{mccalla2013spots,vandenberg2015Rigorous}.  The algebraic matching condition \eqref{MatchEq-cubic} possesses $\mathbb{Z}_2$ symmetry, so $-\mathbf{U}_{R}(r,\theta)$ is also a solution of \eqref{eqn:Galerk} and ring-type patterns emerge in a pitchfork bifurcation. A numerical study of these ring-type patterns first requires solving \eqref{eqn:Ginz} in order to construct an effective initial guess; we will restrict our numerical simulations in Section~\ref{s:numeric} to the simpler spot A-type patterns.  

\subsubsection{Qualitative predictors}
We conclude this section by highlighting the following quantities,
\begin{equation}
    P_1 := c_0, \qquad\qquad P_2 := \frac{[\hat{U}_0]_2}{[\hat{U}_0]_1}, \qquad\qquad P_3 := \frac{[\hat{U}_0]_{1}}{\gamma},  \qquad\qquad P_4 := c_3,
\end{equation}
{  where the notation $[\mathbf{v}]_{i}$ denotes the $i$th element of a vector $\mathbf{v}$.} Each quantity can be computed directly from the original model \eqref{eqn:RD-gen} and describes a qualitative property of our localised solutions. We recall,
\begin{enumerate}[label=(\roman*)]
    \item $P_1$ determines the direction of bifurcation from the Turing point. If $P_1>0$ then localised solutions bifurcate for $\varepsilon>0$ and if $P_1<0$ then localised solutions bifurcate for $\varepsilon<0$.
    \item $P_2$ determines the relationship between $u$ and $v$. If $P_2>0$ then $u$ and $v$ share the same polarity and if $P_2<0$ then $u$ and $v$ have opposing polarity.
    \item $P_3$ determines the polarity of $u$. If $P_3>0$ then spot A-type patterns are made up of peaks and if $P_3<0$ then spot A-type patterns are made up of gaps.
    \item $P_4$ determines whether or not ring-type patterns exist. If $P_4<0$ then ring-type patterns bifurcate from the Turing point, and if $P_4>0$ then no ring-type patterns emerge.
\end{enumerate}
Let us briefly discuss some intuition for each predictor. First, the sign of $P_1$ fixes the sign of $\varepsilon$ via the condition that $\sigma(\lambda;\mu_*+\varepsilon)$ has no real roots for $|\varepsilon|\ll1$. One can formally derive this by a linear transformation of $\sigma(\lambda;\mu_*+\varepsilon)$ into $\{\hat{U}_0,\hat{U}_1\}$ coordinates, so that
\begin{equation*}
    \sigma(\lambda;\mu_*+\varepsilon) = \det\begin{pmatrix}
        \lambda + k^2 + \mathcal{O}(\varepsilon) & -k^2 + \mathcal{O}(\varepsilon) \\ 4 P_1 \varepsilon + \mathcal{O}(\varepsilon^2) & \lambda + k^2 + \mathcal{O}(\varepsilon)
    \end{pmatrix}=0,\quad \implies\quad \lambda = -k^2 \pm 2\mathrm{i}k\sqrt{P_1 \varepsilon} + \mathcal{O}(\varepsilon).
\end{equation*}
Then, we see that $\lambda$ becomes complex when $P_1 \varepsilon>0$. The quantities $P_2, P_3$ both appear in the profile of our spot A type solutions \eqref{RadialProfile}, which we can write as
\begin{equation*}
\mathbf{U}_{A}(r,\theta) = 2\sqrt{3}\; k \; P_3 \;(P_1 \varepsilon)^{\frac{1}{2}}\sum_{n=-N}^{N} a_{|n|}\,J_{|m n|}(k r) \cos(m n \theta) \begin{pmatrix}
    1 \\ P_2
\end{pmatrix} + \mathcal{O}(|\varepsilon|),
\end{equation*}
and so the role of both $P_2, P_3$ can be seen explicitly. The final predictor $P_4$ has already been discussed previously, where we observed that it corresponds to a focussing/defocussing condition in the non-autonomous Ginzburg--Landau equation \eqref{eqn:Ginz}. We will see that these predictors provide a quick guide by which to compare and categorise different reaction--diffusion systems with regards to the formation of localised dihedral patterns.

\subsection{Vegetation models}\label{s:Theory-veg}
We now introduce several models for dryland vegetation and compute the quantities $P_j$ for $j=1,\dots, 4$ in each model. This will then allow us to predict qualitative properties of the types of localised patterns emerging in each model. {  We again recall that our qualitative predictions are regarding the \emph{existence} of localised steady-states, which may represent transient states between uniform and patterned environments. We then support our analysis with a numerical exploration of some spot A-type patterns in Section~\ref{s:numeric}, where we verify our qualitative predictions and observe how localised patterns evolve in each of the following vegetation models.}
\subsubsection{Klausmeier--Gray--Scott model}\label{ss:KG}
We begin with the system
\begin{equation}\label{eqn:KGS}
    \begin{split}
        u_t &= \Delta u + v u^2 - m u,\\
        v_t &= \delta_v \Delta v + \mu - v - v u^2,
    \end{split}
\end{equation}
introduced by van der Stelt et al. in \cite{vanderstelt2013} and commonly referred to as the `extended Klausmeier' model \cite{Siteur2014Beyond,Siero2015}, or the `Klausmeier--Gray--Scott' model \cite{Gandhi2018GrayScott,Wang2021KGS,Li2022Bifurcation,Zelnik2018tristability}. The first name is in reference to the Klausmeier model for vegetation patterns on sloped domains,
\begin{equation*}
    \begin{split}
        u_t &= \Delta u + v u^2 - m u,\\
        v_t &= \nu \partial_x v + \mu - v - v u^2,
    \end{split}
\end{equation*}
where $m>0$ models an effective death rate of vegetation and $\nu$ models the advection of water down a uniform hillslope \cite{Klausmeier1999}. The Klausmeier model can be extended to flat land by setting $\nu=0$ and including a diffusive term for the water density, thus resulting in \eqref{eqn:KGS}. The name `Klausmeier--Gray--Scott' refers to the similarities between \eqref{eqn:KGS} and the Gray--Scott model for chemical reactions \cite{Gray1983Autocatalytic}. We will presently introduce another model derived as an extension of the Klausmeier model, and so we will call \eqref{eqn:KGS} the Klausmeier--Gray--Scott model to reduce confusion. 

\begin{figure}[t!]
    \centering
    \includegraphics[width=\linewidth]{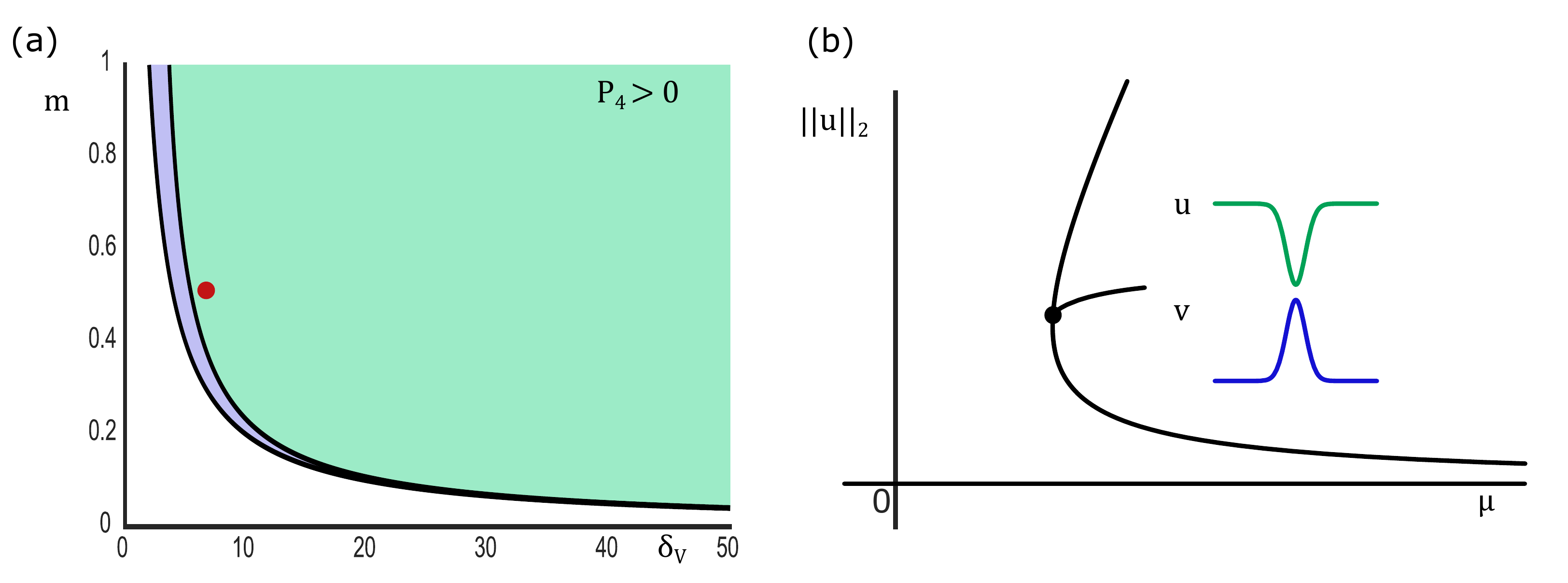}
    \caption{(a) Plot of the sign of $P_4$ in $(\delta_v,m)$-parameter space; $P_4<0$ in the left-most shaded region (blue) and $P_4>0$ in the right-most shaded region (green). The left-most boundary is the curve $\delta_v m = 2$, and so no Turing bifurcation occurs in the unshaded region. The red dot represents the parameter values \eqref{pars:KGS} chosen for our numerical study. (b) Schematic bifurcation diagram of the Klausmeier--Gray--Scott model \eqref{eqn:KGS} with parameter values \eqref{pars:KGS}, where anti-phase spot A-type patterns emerge as gaps from the Turing point as $\mu$ increases.}
    \label{fig:K-GS}
\end{figure}

The Klausmeier--Gray--Scott model \eqref{eqn:KGS} is described by the general reaction--diffusion model \eqref{eqn:RD-gen} with $D_v = \delta_v$, $\beta=0$, and
\begin{equation*}
    \begin{split}
        \hat{f}(u,v;\mu) &= - v u^2 + m u,\qquad
        \hat{g}(u,v;\mu) = - \mu + v + v u^2.
    \end{split}
\end{equation*}
For this particular model, the functions $\hat{f},\hat{g}$ are simple enough to analyse explicitly, which is done in Appendix~\ref{app:KGS-analysis}. We discover that \eqref{eqn:KGS} undergoes a Turing bifurcation if $\delta_v m >2$, and 
\begin{equation*}
    P_1 > 0, \qquad\qquad P_2<0, \qquad\qquad P_3<0 \qquad\qquad \text{for all }\; m >0,\quad \delta_v > 2/m.
\end{equation*}
The sign of $P_4$ varies for different values of $\delta_v, m$, and so we numerically plot a diagram for $\text{sgn}\,(P_4)$ in  Figure~\ref{fig:K-GS}(a). 

The benefit of finding explicit expressions for each $P_j$ is that we can judge the robustness of each qualitative property; i.e., how sensitive any numerical solutions are with respect to our choice of parameters. What we find is that the direction of bifurcation ($P_1$), whether solutions are in-phase or anti-phase ($P_2$), and whether the solutions are peaks or gaps ($P_3$) remain consistent across all choices of $\delta_v,m$ such that $\delta_v m>2$. Figure~\ref{fig:K-GS}(a) shows that the subcriticality condition $P_4<0$ only holds for sufficiently small values of $\delta_v$, given a fixed $m$, and so the ring-type patterns found in Theorem~\ref{thm:Ring} emerge in a much smaller parameter region than the spot A-type patterns of Theorem~\ref{thm:Patch}. 

For our subsequent numerical study in Section~\ref{s:numeric}, we choose the following parameter values
\begin{equation}\label{pars:KGS}
    m = 0.5,\quad \delta_v = 7.2,
\end{equation}
taken from \cite{Zelnik2018tristability}, where localised periodic states have been observed previously; we will use these parameter values for the remainder of this work when considering \eqref{eqn:KGS}. Then, we identify the Turing point
\begin{equation*}
    (u_*,v_*,\mu_*) = (1.071,\; 0.467,\; 1.002),
\end{equation*}
with
\begin{equation*}
    P_1 = 6.923, \qquad\qquad P_2 = -0.348, \qquad \qquad P_3 = -1.512, \qquad\qquad P_4 = 0.248,
\end{equation*}
which is consistent with our explicit calculations. Hence, the Klausmeier--Gray--Scott model \eqref{eqn:KGS} with parameter values \eqref{pars:KGS} exhibits anti-phase spot A-type localised gaps, which bifurcate as the precipitation increases beyond the Turing point; see Figure~\ref{fig:K-GS}(b).

We note that the direction of bifurcation, as predicted by $P_1$, does not tell the whole story; the bifurcation curve in Figure~\ref{fig:K-GS}(b) could undergo a fold and continue into the region where $\mu<\mu_*$. This is the case for certain localised one dimensional solutions in the Klausmeier--Gray--Scott model \cite{Zelnik2018tristability}, and so we might expect the same to be true in this case.

\subsubsection{Logistic Klausmeier model}\label{ss:Kl}
Our second example is the logistic Klausmeier model introduced by Bastiaansen et al. \cite{Bastiaansen2019Stable}, 
\begin{equation}\label{eqn:extK}
    \begin{split}
        u_t &= \Delta u + (1-b u) v u^2 - m u,\\
        v_t &= \delta_v \Delta v + \mu - v - v u^2,
    \end{split}
\end{equation}
which is an extension of the previous Klausmeier--Gray--Scott model \eqref{eqn:KGS}, where the unrestricted vegetation growth term $v u^2$ is replaced by a logistic growth term $(1 - b u) v u^2$ with carrying capacity $b^{-1}>0$. 

This adjustment limits the growth of vegetation at any given point, which is motivated as follows. For extremely low precipitation levels it is reasonable to assume that vegetation growth is constrained by the amount of soil water, as modelled by \eqref{eqn:KGS}; however, in environments that can support larger quantities of vegetation, one should expect other limiting effects to become relevant, such as competition for space and nutrients among individual plants. Hence, when considering larger scale vegetation patterns, such as in \cite{Bastiaansen2019Stable,byrnes2023large,Carter2018,Sewalt2017Multipulse,Carter2023instability,Iuorio2021autotoxicity}, it may be more appropriate to consider \eqref{eqn:extK} rather than \eqref{eqn:KGS}. 

\begin{figure}[t!]
    \centering
    \includegraphics[width=\linewidth]{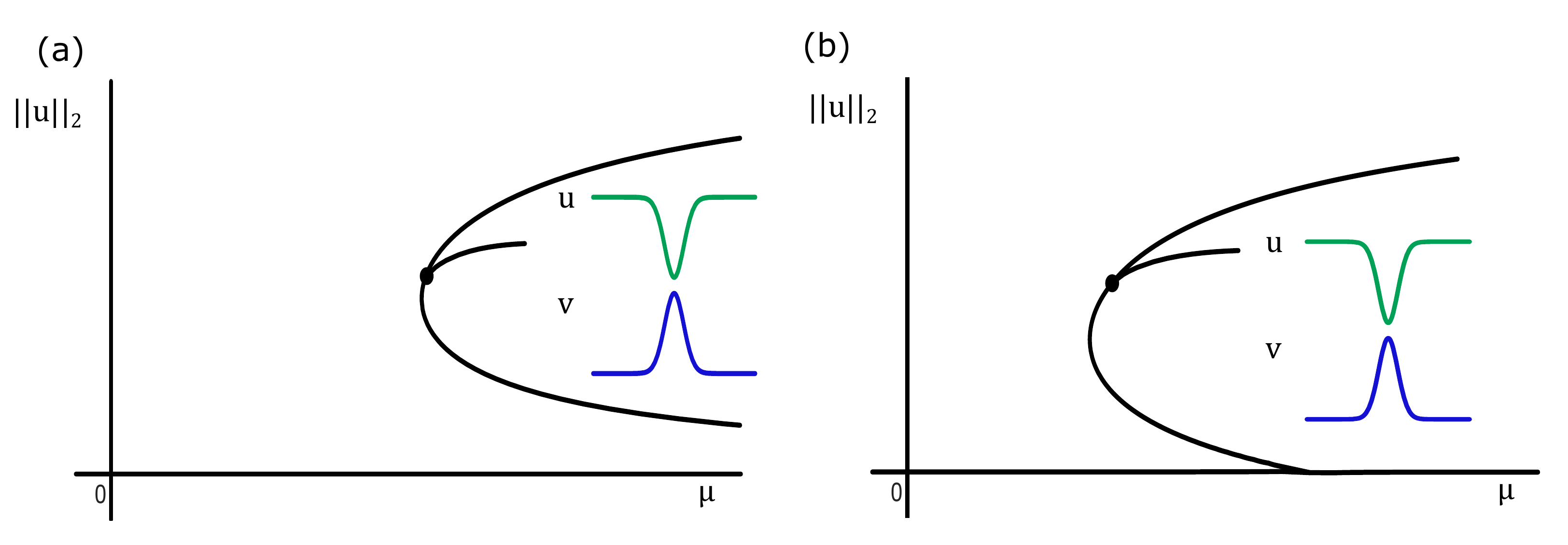}
    \caption{Schematic bifurcation diagrams of (a) the logistic Klausmeier model \eqref{eqn:extK} with parameter values \eqref{pars:extK}, and (b) the NFC--Gilad model \eqref{eqn:Gilad} with parameter values \eqref{pars:Gilad}. In both cases, anti-phase spot A-type patterns emerge as gaps from the Turing point as $\mu$ increases.}
    \label{fig:Kl-Gi}
\end{figure}

The logistic Klausmeier model \eqref{eqn:extK} is again covered by the general reaction--diffusion model \eqref{eqn:RD-gen} with $D_v = \delta_v$, $\beta=0$, and
\begin{equation}
    \begin{split}
        \hat{f}(u,v;\mu) &= -(1-b u) v u^2 + m u,\qquad
        \hat{g}(u,v;\mu) = - \mu + v + v u^2.
    \end{split}
\end{equation}
 The additional nonlinear term makes any explicit calculations significantly more cumbersome, and so we just present our numerical calculations of $P_j$. We fix $b=1$ and choose the following parameter values
\begin{equation}\label{pars:extK}
    m = 0.45,\quad \delta_v = 182.5,
\end{equation}
used to model grass patterns in \cite{Klausmeier1999}, where the value of the diffusion constant $\delta_v$ is significantly larger than in the previous example. Then, we find the Turing point
\begin{equation*}
    (u_*,v_*,\mu_*) = (0.465,\; 1.809,\; 2.200),
\end{equation*}
with
\begin{equation*}
    P_1 = 0.503, \qquad\qquad P_2 = -0.282, \qquad \qquad P_3 = -0.965, \qquad\qquad P_4 = 0.015.
\end{equation*}
Hence, the logistic Klausmeier model \eqref{eqn:extK} with parameter values \eqref{pars:extK} exhibits anti-phase spot A-type localised gaps, which bifurcate as the precipitation increases beyond the Turing point. Comparing Figure~\ref{fig:Kl-Gi}(a) with Figure~\ref{fig:K-GS}(b), we see that the localised structures emerging in the Klausmeier--Gray--Scott model \eqref{eqn:KGS} and the logistic Klausmeier model \eqref{eqn:extK} are qualitatively the same.

\subsubsection{NFC--Gilad model}\label{ss:Gi}
Our next model is a simplification of the Gilad model introduced in \cite{Gilad2004Engineers,Gilad2007}
\begin{equation*}
    \begin{split}
        u_t ={}& G_u u(1-u) - u + \delta_u\Delta u,\\
        v_t ={}& \alpha\left(\frac{u + \gamma \sigma}{u + \gamma}\right)w - { \nu}(1 - \rho u)v - G_v v + \delta_v \Delta v,\\
        w_t ={}& \mu - \alpha\left(\frac{u + \gamma \sigma}{u + \gamma}\right) w + \delta_w \Delta(w^2),\\
    \end{split}
\end{equation*}
where $w(t,\mathbf{x})$ denotes the density of water on the surface of the soil and $G_u(t,\mathbf{x}), G_v(t,\mathbf{x})$ are non-local terms that describe vegetation growth and soil water consumption, {  given by 
\begin{equation*}
    \begin{split}
        G_{u}(t,\mathbf{x}) ={}& \frac{\Lambda}{2\pi}\int_{\Omega} \frac{1}{S_0^2}\exp\left(-\frac{|\mathbf{x}-\mathbf{x}'|}{2 S_0^2 (1 + \eta u(t,\mathbf{x}))^2}\right)\,v(t,\mathbf{x}')\,\mathrm{d}\mathbf{x}',\\
        G_{v}(t,\mathbf{x}) ={}& \frac{\Lambda}{2\pi}\int_{\Omega} \frac{1}{S_0^2}\exp\left(-\frac{|\mathbf{x}-\mathbf{x}'|}{2 S_0^2 (1 + \eta u(t,\mathbf{x}))^2}\right)\,u(t,\mathbf{x}')\,\mathrm{d}\mathbf{x}',\\
    \end{split}
\end{equation*}
respectively}. One of the aims of this model was to describe the behaviour of plants with widely spread out roots, such that the interaction between vegetation and water is no longer strictly local, resulting in a three-component non-local PDE system. As such, this model is not covered by our results unless some simplifications are first performed.  

We consider a simplification of the Gilad model introduced by Zelnik et al. in \cite{Zelnik2015Gradual}, where the authors focused on studying fairy circles on flat terrains in Namibia. For this particular environment, the sandy soil results in very fast changes to the surface water $w(t,\mathbf{x})$---which converges to an equilibrium state---and the remaining slow dynamics of the system can be reduced to the two-component reaction-diffusion system
\begin{equation}\label{eqn:Gilad}
\begin{split}
    u_t &= \Delta u + \Lambda v u (1 - u)(1+\eta u)^2 - u,\\
    v_t &= \delta_{v}\Delta v + \mu - \nu(1- \rho u)v - \Lambda v u (1+\eta u)^2.\\
\end{split}\end{equation}
We refer to this model as the `Namibian fairy circle Gilad' model, or `NFC--Gilad' model for short. It is exactly of the form in \eqref{eqn:RD-gen}, with $D_v=\delta_v$, $\beta=0$, and
\begin{equation}\begin{split}
    \hat{f}(u,v;\mu) &= - \Lambda v u (1 - u)(1+\eta u)^2 + u,\qquad
    \hat{g}(u,v;\mu) = - \mu + \nu(1- \rho u)v + \Lambda v u (1+\eta u)^2.\\
\end{split}\end{equation}
Again, the functions $\hat{f},\hat{g}$ are complicated enough that we only present numerical calculations here. We choose the following parameter values
\begin{equation}\label{pars:Gilad}
    \Lambda  = \tfrac{16}{35}, \quad \eta = \tfrac{14}{5}, \quad \nu = \tfrac{10}{7}, \quad \rho = \tfrac{7}{10}, \quad \delta_v = 125,
\end{equation}
taken from \cite{Zelnik2015Gradual}; then, we find the Turing point
\begin{equation*}
    (u_*,v_*,\mu_*) = (0.474,\; 0.768,\; 1.635),
\end{equation*}
with qualitative predictors
\begin{equation*}
    P_1 = 0.381, \qquad\qquad P_2 = -0.207, \qquad \qquad P_3 = -0.575, \qquad\qquad P_4 = 0.818.
\end{equation*}
Hence, the NFC--Gilad model \eqref{eqn:Gilad} with parameter values \eqref{pars:Gilad} exhibits anti-phase spot A-type localised gaps, which bifurcate as the precipitation increases beyond the Turing point; see Figure~\ref{fig:Kl-Gi}(b). Again, comparing Figure~\ref{fig:Kl-Gi}(b) with Figures~\ref{fig:K-GS}(b) \& \ref{fig:Kl-Gi}(a), we see that the localised structures emerging in the NFC--Gilad model \eqref{eqn:Gilad} are qualitatively the same as those in the Klausmeier--Gray--Scott and logistic Klausmeier models.

\subsubsection{von Hardenberg model}\label{ss:vH}
The final model we consider is the von Hardenberg model introduced in \cite{vonHardenberg2001},
\begin{equation}\label{eqn:vonH}
    \begin{split}
        u_t &= \Delta u + \tfrac{\gamma v}{1 + \sigma v}u - u^2 - \nu u,\\
        v_t &= \delta_{v}\Delta (v-\beta u) + \mu - (1- \rho u)v - u v^2,
    \end{split}
\end{equation}
which reproduces many symmetry-breaking transitions found in field observations. In particular, simulations of domain-covering patterns undergo `gaps $\to$ labyrinth $\to$ spots' and `spots $\to$ stripes $\to$ gaps' transitions \cite{Gowda2014,Gowda2016}, and simulations of localised patterns undergo `spot $\to$ ring' and `ring $\to$ spots' transitions \cite{Meron2004}.  

We highlight some of the terms in \eqref{eqn:vonH} not seen in the previous models. The {  coefficient $\left(\frac{\gamma v}{1 + \sigma v}\right)$} models vegetation growth due to water intake, which is linear for low levels of water and tends to a constant for larger levels of water; this is instead of the constant term { $m$} in \eqref{eqn:KGS} and \eqref{eqn:extK} or the linear term { $\Lambda v$} in \eqref{eqn:Gilad}. The $\Delta(v - \beta u)$ term models the transport of water via Darcy's law; this includes the suction of water by plant roots at a rate of $\beta>0$, which is not covered by the previous models. 
 
\begin{figure}[t!]
    \centering
    \includegraphics[width=\linewidth]{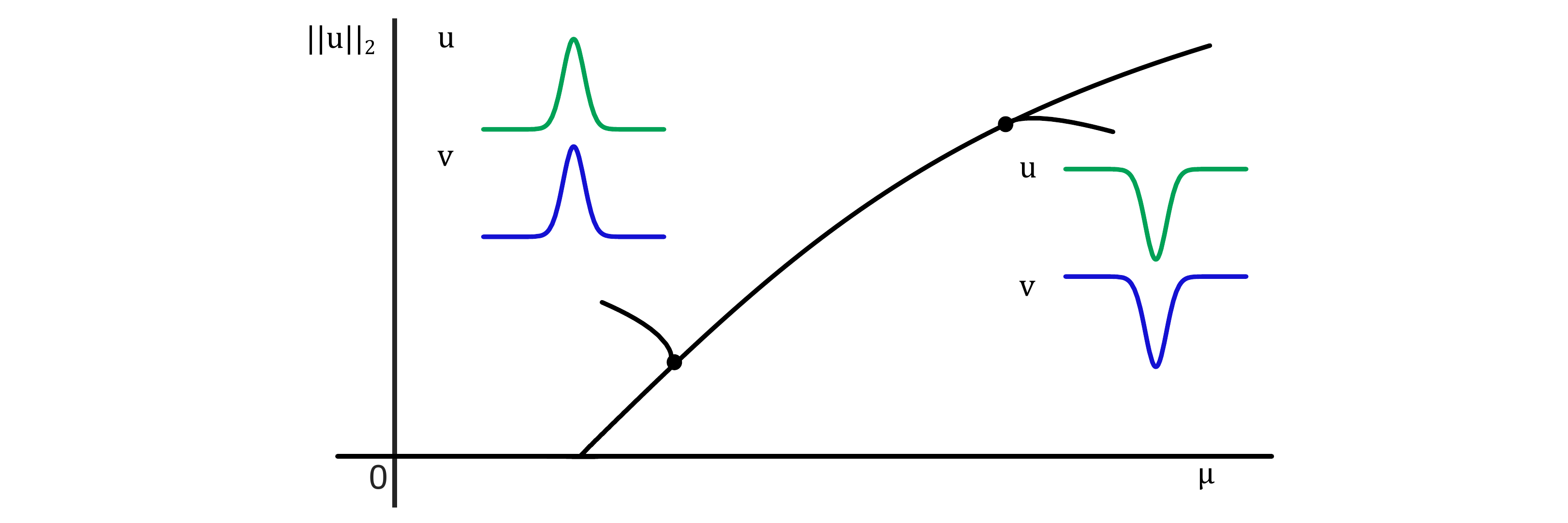}
    \caption{Schematic bifurcation diagram of the von Hardenberg model \eqref{eqn:vonH} with parameter values \eqref{pars:vonH}, with in-phase spot A-type patterns emerging as peaks from the lower Turing point \eqref{Turing:vH-1} and as gaps from the upper Turing point \eqref{Turing:vH-2}.}
    \label{fig:v-H}
\end{figure}

The von Hardenberg model \eqref{eqn:vonH} is of the form in \eqref{eqn:RD-gen}, with $D_v=\delta_v$, $\beta>0$, and
\begin{equation}
    \begin{split}
        \hat{f}(u,v;\mu) &= - \tfrac{\gamma v}{1 + \sigma v}u + u^2 + \nu u,\qquad
        \hat{g}(u,v;\mu) = - \mu + (1- \rho u)v + u v^2,
    \end{split}
\end{equation}
and we again only present numerical calculations for $P_1, P_2, P_3, P_4$ here. We choose the following parameter values
\begin{equation}\label{pars:vonH}
    \gamma = 1.6,\quad \sigma = 1.6,\quad \nu = 0.2,\quad \rho = 1.5,\quad\delta_v = 100,\quad \beta = 3,
\end{equation}
taken from \cite{vonHardenberg2001}. Notably, we obtain two Turing points for the von Hardenberg model, both with qualitatively different properties. The first Turing point is at
\begin{equation}\label{Turing:vH-1}
    (u_*,v_*, \mu_*) = (0.017, \; 0.173, \; 0.169),
\end{equation}
with qualitative predictors
\begin{equation*}
    P_1 = -0.384, \qquad\qquad P_2 = 1.707, \qquad \qquad P_3 = 0.8427, \qquad\qquad P_4 = 0.0012,
\end{equation*}
while the second Turing point is at
\begin{equation}\label{Turing:vH-2}
    (u_*,v_*, \mu_*) = (0.271, \; 0.556, \; 0.414),
\end{equation}
with qualitative predictors
\begin{equation*}
    P_1 = 0.217, \qquad\qquad P_2 = 2.578, \qquad \qquad P_3 = -1.512, \qquad\qquad P_4 = 0.014.
\end{equation*}

Hence, the von Hardenberg model \eqref{eqn:vonH} with parameter values \eqref{pars:vonH} exhibits the following localised dihedral patterns. As the precipitation decreases beyond the first Turing point \eqref{Turing:vH-1}, in-phase spot A-type localised peaks emerge. In contrast, as the precipitation increases beyond the second Turing point \eqref{Turing:vH-2}, in-phase spot A-type localised gaps emerge; for a summary of these predictions, see Figure~\ref{fig:v-H}. We note that, while ring-type patterns do not emerge from either Turing point for the given parameter values, localised one dimensional patterns have been found to emerge from the first Turing point for a different choice of parameters; see \cite[Figure~15]{Dawes2016}. As discussed earlier, these one dimensional patterns require $c_3<0$ in order to emerge, and so we expect dihedral ring-type patterns to also emerge in those parameter regimes.

\section{Numerical Results}\label{s:numeric}

In the previous section we demonstrated how fully localised patterns emerge in general two-component reaction--diffusion systems near a Turing instability, {  which we now support with a brief numerical study on the evolution of spot A-type dihedral vegetation patterns in each of the vegetation models introduced previously. The goal of this section is to numerically explore spot A-type vegetation patterns in the four vegetation models introduced previously, in order to (1) verify our qualitative predictions from Section \ref{s:Theory-veg}; (2) explore the role of localised patterns as transient states of vegetation models; and (3) observe how different models with similar initial conditions can result in qualitatively different evolution dynamics.

The numerical codes for producing the results of this section are all available at \cite{Hill2023Github-Veg}, along with codes for calculating the predictive quantities $P_1,\dots,P_4$ in each vegetation model considered here. Furthermore, we provide others with the option to define their own reaction-diffusion models of the form \eqref{eqn:RD-gen} and determine the qualitative properties of patterns emerging from a given Turing point.}

\subsection{Implementation}
Our numerical approach is as follows: for a system of the form \eqref{eqn:RD-gen} we provide the functions $\hat{f}(u,v,\mu)$, $\hat{g}(u,v,\mu)$, including the values of any other parameters, and an initial guess for the Turing point $(u_*, v_*, \mu_*)$. We numerically solve the algebraic equations detailed in Section~\ref{s:Theory} in order to find a Turing point close to our initial guess, and compute $k, \hat{U}_0, \hat{U}_{1}, \gamma, c_0, c_3$. Then the the values of $P_1, P_2, P_3, P_4$ can be computed, as in the previous section for each vegetation model.  

Next, we choose what type of dihedral pattern we want to find by fixing $m,N\in\mathbb{N}$ in Theorem~\ref{thm:Patch} and providing an initial guess for the algebraic matching condition \eqref{MatchEq} for spot A-type patterns. This results in an initial profile \eqref{RadialProfile} for a localised dihedral pattern bifurcating from the Turing instability; we then input this solution into our time--stepping codes in order to investigate how these solutions evolve over time. 

We employ a second-order exponential time--differencing method from Asant\'e--Asamani et al.~\cite{Asante-Asamani2020ETD} for a reaction--diffusion system of the form
\begin{equation}
    \partial_t \mathbf{u} = \begin{pmatrix}
        D_u & 0  \\ 0 & D_v
    \end{pmatrix}\Delta\mathbf{u} + \mathbf{F}(\mathbf{u};\mu),
\end{equation}
with homogeneous Neumann boundary conditions. In order to convert our general system \eqref{eqn:RD-gen} into this form, we apply the transformation
\begin{equation}\label{transf}
        \hat{\mathbf{u}}:=\begin{pmatrix}
        1 & 0\\ -\tfrac{\beta D_v}{D_v -1} & 1
    \end{pmatrix}\mathbf{u}, \qquad\qquad  \mathbf{u}=\begin{pmatrix}
        1 & 0\\ \tfrac{\beta D_v}{D_v -1} & 1
    \end{pmatrix}\hat{\mathbf{u}},
\end{equation}
so that our general system of the form
\begin{equation}
    \partial_t \mathbf{u} = \begin{pmatrix}
        1 & 0  \\ -D_v \beta & D_v
    \end{pmatrix}\Delta\mathbf{u} + \mathbf{F}(\mathbf{u};\mu)
\end{equation}
becomes
\begin{equation}
    \partial_t \hat{\mathbf{u}} = \begin{pmatrix}
        1 & 0  \\ 0 & D_v
    \end{pmatrix}\Delta\hat{\mathbf{u}} + \begin{pmatrix}
        1 & 0\\ -\tfrac{\beta D_v}{D_v -1} & 1
    \end{pmatrix}\mathbf{F}\left(\begin{pmatrix}
        1 & 0\\ \tfrac{\beta D_v}{D_v -1} & 1
    \end{pmatrix}\hat{\mathbf{u}};\mu\right).
\end{equation}
In particular, we apply our time-stepping codes to the system
\begin{equation}
    \partial_t \begin{pmatrix} \hat{u} \\ \hat{v}\end{pmatrix} = \begin{pmatrix}
       1 & 0  \\ 0 & D_v
    \end{pmatrix}\Delta\begin{pmatrix} \hat{u} \\ \hat{v}\end{pmatrix} + \begin{pmatrix} -\hat{f}\left(\hat{u},\hat{v} + \tfrac{\beta D_v \hat{u}}{D_v-1};\mu\right) \\ \left(\tfrac{\beta D_v}{D_v -1}\right)\hat{f}\left(\hat{u},\hat{v} + \left(\tfrac{\beta D_v}{D_v -1}\right)\hat{u};\mu\right)-\hat{g}\left(\hat{u},\hat{v} + \left(\tfrac{\beta D_v}{D_v -1}\right)\hat{u};\mu\right)\end{pmatrix},
\end{equation}
and then recover $u,v$ from the transformation $u = \hat{u}$, $v = \frac{\beta D_v}{D_v - 1}\hat{u} + \hat{v}$.  

We compute our simulations on a square domain of length $20\lambda$, where $\lambda:=\frac{2\pi}{k}$ is the wavelength related to each Turing bifurcation. This $k-$dependent choice of domain means that any patterns found across different models should be of a similar size when plotted. We present our simulations through snapshots of our solutions in time; videos of each simulation can be found at \cite{Hill2023Github-Veg}.
\subsection{Simulations}
\begin{figure}[t!]
    \centering
    \includegraphics[width=\linewidth]{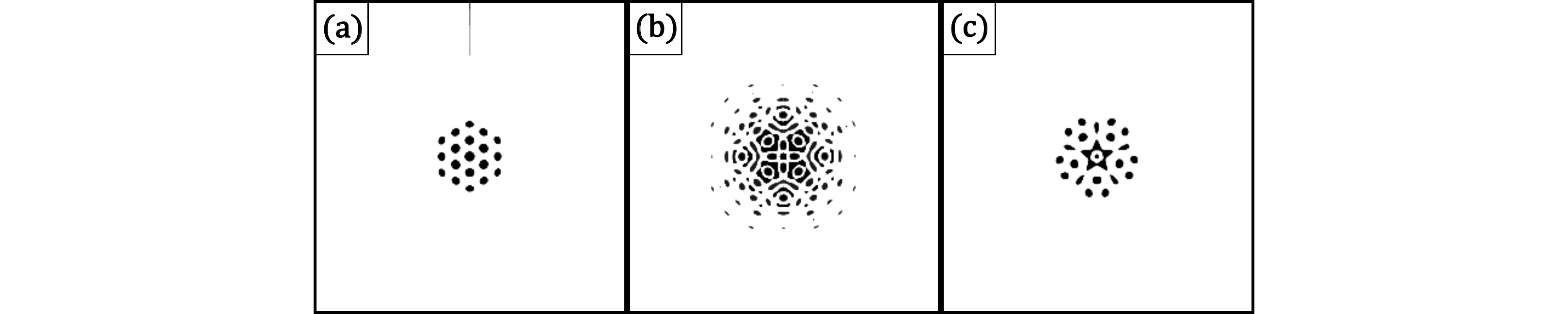}
    \caption{Three localised dihedral patterns: (a) hexagons $(\mathbb{D}_{6})$, (b) squares $(\mathbb{D}_{4})$, and (c) pentagons $(\mathbb{D}_{5})$. Each profile is an approximate localised dihedral pattern obtained in \cite{Hill2023DihedralPatch} and together they form the initial guesses for our simulations.}
    \label{fig:init}
\end{figure} 
We present three examples of localised dihedral patterns for our numerical study: hexagons ($\mathbb{D}_{6}$), squares ($\mathbb{D}_{4}$) and pentagons ($\mathbb{D}_{5}$). In particular, we provide an initial guess
\begin{equation}\label{InitialGuess}
\begin{pmatrix}
    u\\v
\end{pmatrix} =\begin{pmatrix}
    u_*\\v_*
\end{pmatrix} +  C P_3 \left[\sum_{n=-N}^{N} a_{|n|}\,J_{|m n|}(k r) \cos(m n \theta)\right]\mathrm{e}^{-(P_1 \varepsilon)^{\frac{1}{2}} r}\begin{pmatrix}
    1 \\ P_2
\end{pmatrix}
\end{equation} 
adapted from \eqref{RadialProfile}, where $\varepsilon = \mu - \mu_*$, $C>0$ is a fixed constant, and 
\begin{equation*}\begin{split}
    \text{Hexagon:} & \qquad m=6, \qquad N=2, \qquad (a_0,a_1,a_2) = (0.311, \; 0.267, \; 0.189),\\
    \text{Square:} & \qquad m=4, \qquad N=5, \qquad (a_0,a_1,a_2,a_3,a_4,a_5) = (-0.136, \; 0.262, \; 0.236, \; -0.114, \; 0.187, \; 0.145),\\
    \text{Pentagon:} & \qquad m=5, \qquad N=3, \qquad (a_0,a_1,a_2,a_3) = (-0.382, \; 0.300, \; 0.382,\; 0.486).
\end{split}\end{equation*}
See Figure~\ref{fig:init} for the profile of $u-u_*$ in each case. Note that, while we refer to each pattern by the polygon associated with its dihedral symmetry group---i.e., hexagons, squares and pentagons---the patterns do not necessarily appear similar to these polygons. For example, the `square' pattern is related to a twelve-fold quasipattern, rather than a square lattice.  

These three examples of dihedral patterns each possess a different level of rarity in pattern formation; hexagons are very common in experiments and field observations, square lattices are common but our specific square example is not, and pentagons are not very common at all. However, Theorem~\ref{thm:Patch} tells us that each example emerges from the same Turing point {  in the same parameter regime, and so it would be interesting to compare the time evolution for each example and observe whether the square and pentagon patterns converge to the more common hexagons over time.}

{  Since each of the models we consider is nondimensionalised, we cannot readily compare the amplitude of patterns across different models. However, we observe that, for a given model, each of our initial dihedral patterns converge to roughly the same amplitude, and so we provide an approximate range of values for $u$ in each model. The videos available at \cite{Hill2023Github-Veg} show an individual colorbar for each pattern, thus providing more details for enquiring readers.}

\subsubsection{Klausmeier--Gray--Scott}
\begin{figure}[t!]
    \centering
    \includegraphics[width=0.9\linewidth]{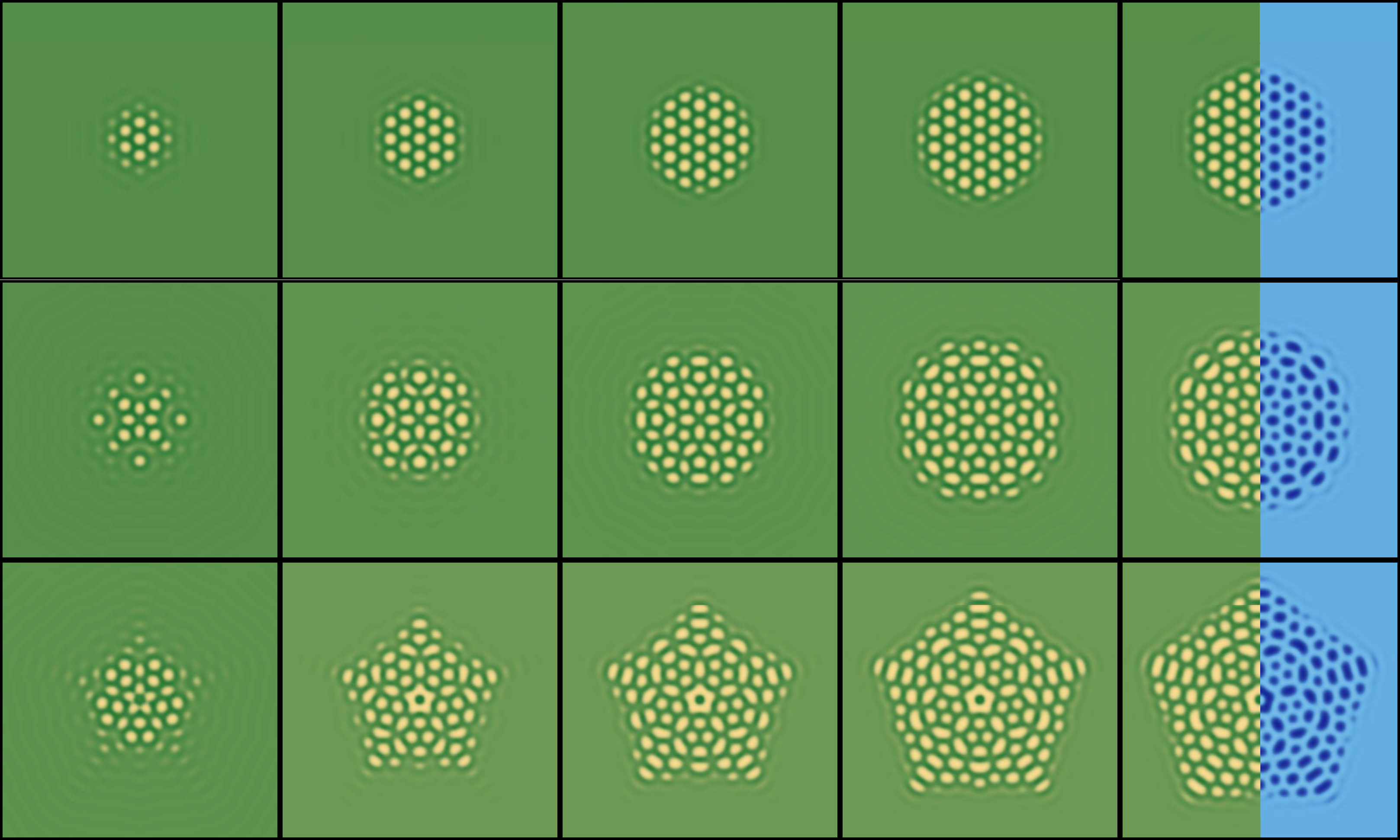}
    \caption{Time evolution of vegetation density $u$ in the Klausmeier--Gray--Scott model \eqref{eqn:KGS} with parameter values \eqref{pars:KGS}. The dark green represents a vegetated state ({  $u\approx1.5$}), while the lighter yellow represents a bare state ({  $u\approx0$}). Each row corresponds to a different initial perturbation \eqref{InitialGuess} for---from top to bottom---hexagons, squares, and pentagons. From left to right, the columns correspond to different points in time $t=100, 200, 300, 400$ and $500$. {  The final column is divided into vegetation density (left) and water density (right), where darker and lighter blue represents higher and lower water density, respectively. Videos of each simulation can be found at \cite{Hill2023Github-Veg}.}}
    \label{fig:sim-KG}
\end{figure}
We begin our simulations with the Klausmeier--Gray--Scott model \eqref{eqn:KGS} with parameter values \eqref{pars:KGS}. {  We recall from Section \ref{ss:KG} that
\begin{equation*}
    P_1 = 6.923>0,\qquad P_2 = -0.348<0, \qquad P_3 = -1.512<0, \qquad P_4 = 0.248>0
\end{equation*}
and so we expect localised anti-phase gaps to emerge as the precipitation increases from the Turing point $\mu_*=1.002$. Hence, we simulate the time evolution of our initial guess \eqref{InitialGuess}---for hexagons, squares and pentagons---in the parameter region $\mu>\mu_*$.}

The Turing bifurcation has an associated wave number $k\approx 0.3177$, and so we simulate our patterns in a square box of length $396$. The simulations of our three examples---hexagons, square, and pentagons---are presented in Figure~\ref{fig:sim-KG}. We note that the patterns found for \eqref{eqn:KGS} are all anti-phase, as predicted in Section~\ref{s:Theory} {  and shown in the final column of Figure~\ref{fig:sim-KG}}; we do not {  otherwise} plot the water density $v$ since the spatial structure is identical to the vegetation density $u$. 

For the hexagon example, solutions form a regular hexagonal lattice of gaps. As time evolves, new gaps emerge on the outer layer of the patch such that the width of the pattern grows continuously. This behaviour is what we would expect for localised patterns that undergo homoclinic snaking (for example \cite{lloyd2008localized,Hill2023DihedralPatch,Bramburger2019Rolls} for snaking solutions and \cite{Lloyd2019InvasionFronts,Lloyd2021HexInvasionFronts} for growing patterns away from the snaking region) and is also observed in the square and pentagon examples, but with very different structures being formed.  

In the square example, the initial pattern contains regions of the uniform state between individual gaps; emerging gaps first fill these regions until the resulting pattern is compact, before then emerging on the outer layers. Whereas the hexagon example only consists of circular gaps, the square pattern is also comprised of long thin gaps. For later times, the square pattern is reminiscent of the penta-hepta defects and grain boundaries studied in \cite{Subramanian2021Grain}; this is when a point in an almost-hexagonal lattice has five or seven neighbours, rather than the standard six neighbours. 

For the pentagon example, the story is much the same as the square example. Gaps first fill out the initial patterned region before then emerging on the outer layer. Notably, the five gaps at the core of the pattern coalesce and appear to converge to a ring, while the outer layer of the pattern traces out a regular pentagon as times increases.

\subsubsection{Logistic Klausmeier}
\begin{figure}[t!]
    \centering
    \includegraphics[width=0.9\linewidth]{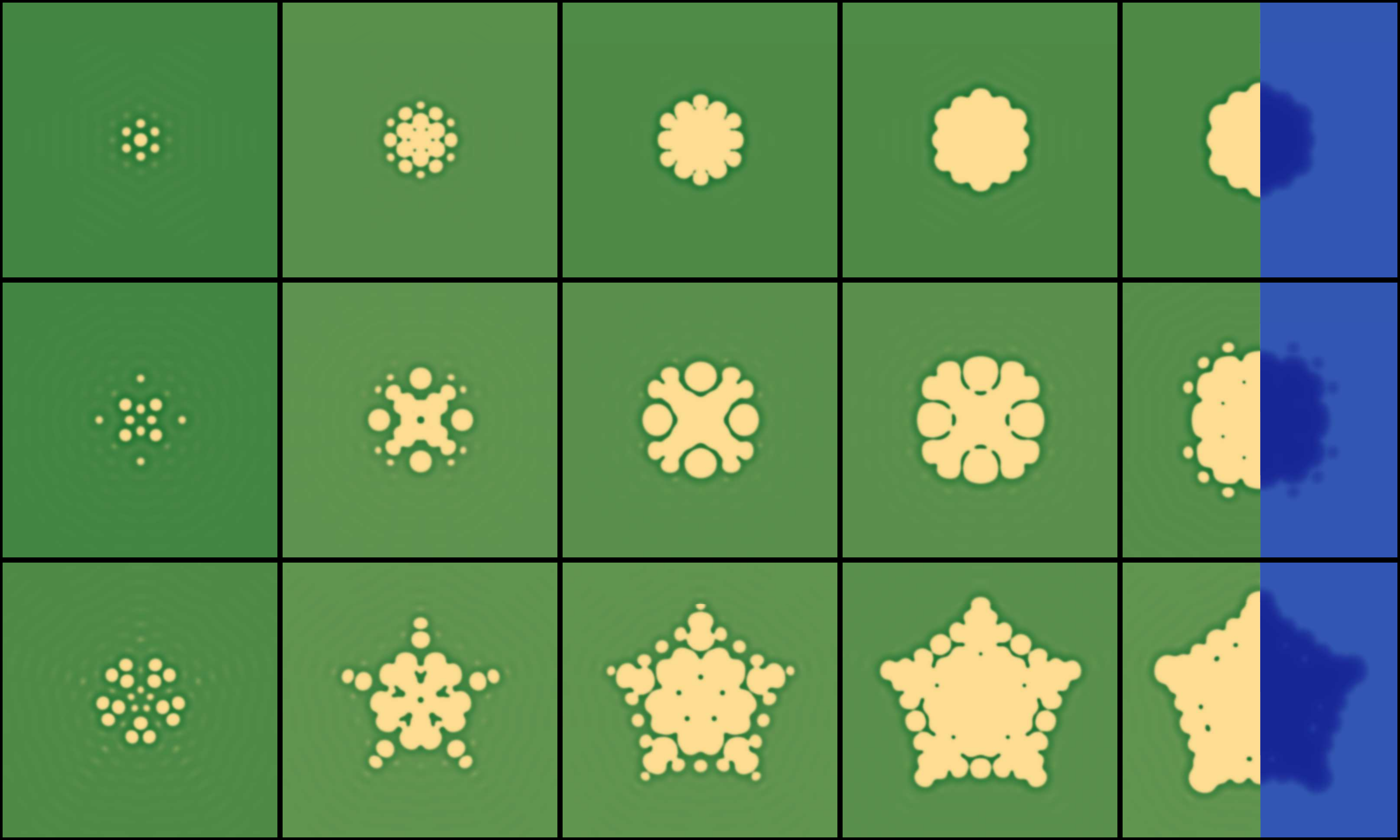}
    \caption{Time evolution of vegetation density $u$ in the logistic Klausmeier model \eqref{eqn:extK} with parameter values \eqref{pars:extK}. The dark green represents a vegetated state ({  $u\approx0.6$}), while the lighter yellow represents a bare state ({  $u\approx0$}). The rows and columns are identical to Figure~\ref{fig:sim-KG}. {  The final column is divided into vegetation density (left) and water density (right), where darker and lighter blue represents higher and lower water density, respectively. Videos of each simulation can be found at \cite{Hill2023Github-Veg}.}}
    \label{fig:sim-Kl}
\end{figure}
{ We next consider} the logistic Klausmeier model \eqref{eqn:extK} with parameter values \eqref{pars:extK}. {  We recall from Section \ref{ss:Kl} that
\begin{equation*}
    P_1 = 0.503 > 0,\qquad P_2 = -0.282<0, \qquad P_3 = -0.965<0, \qquad P_4 = 0.015>0
\end{equation*}
and so we expect localised anti-phase gaps to emerge as the precipitation increases from the Turing point $\mu_*=2.200$. Hence, we simulate the time evolution of our initial guess \eqref{InitialGuess}---for hexagons, squares and pentagons---in the parameter region $\mu>\mu_*$.}

The Turing bifurcation has an associated wave number $k\approx 0.1612$, and so we simulate our patterns in a square box of length $780$. The simulations are presented in Figure~\ref{fig:sim-Kl}, where we again {  mostly} plot the vegetation density $u$ since the water density $v$ exhibits the same spatial structures in anti-phase. 

The logistic term in \eqref{eqn:extK} restricts the height of individual gaps, causing the gaps to grow in width instead; this was also observed in \cite{Bastiaansen2019Stable} for one dimensional patterns. As a result, nearby gaps expand and coalesce, causing a collapse from a patterned state to a compact region of bare soil. New gaps emerge on the periphery of this compact region, in the same positions as in the Klausmeier--Gray--Scott model \eqref{eqn:KGS}, but are soon absorbed into the central gap as time evolves. 

Over time these patterns resemble radial fronts connecting the bare state to the vegetated state, similar to the axisymmetric and one dimensional fronts studied in \cite{Bastiaansen2019Stable,byrnes2023large,Carter2023instability}. However, the geometry of the front interface is highly nontrivial, and so it is unclear whether these structure could be analysed using similar techniques.
\,
\subsubsection{NFC--Gilad}
\begin{figure}[t!]
    \centering
    \includegraphics[width=0.9\linewidth]{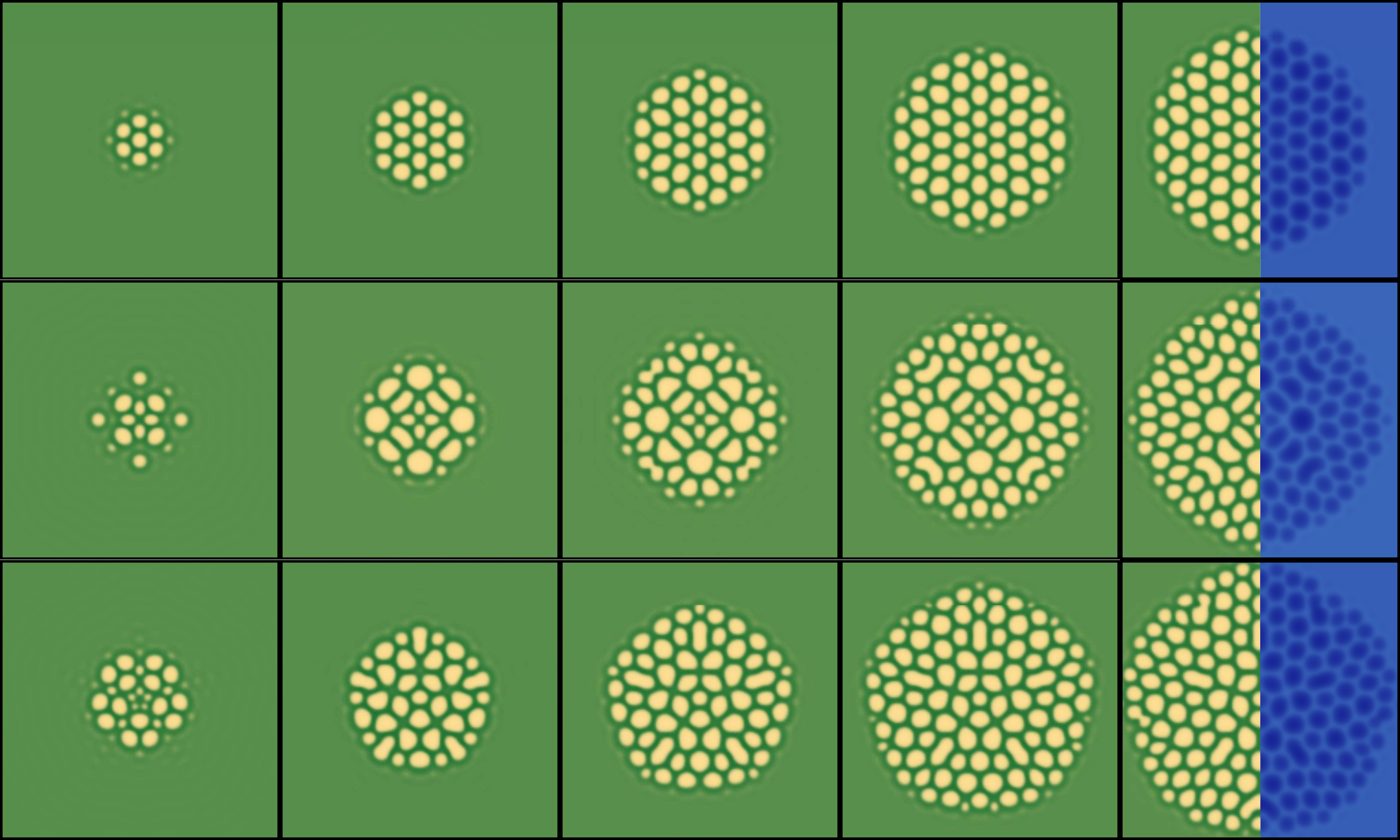}
    \caption{Time evolution of vegetation density $u$ in the NFC--Gilad model \eqref{eqn:Gilad} with parameter values \eqref{pars:Gilad}. The dark green represents a vegetated state ({  $u\approx0.6$}), while the lighter yellow represents a bare state ({  $u\approx0$}). The rows and columns are identical to Figure~\ref{fig:sim-KG}. {  The final column is divided into vegetation density (left) and water density (right), where darker and lighter blue represents higher and lower water density, respectively. Videos of each simulation can be found at \cite{Hill2023Github-Veg}.}}
    \label{fig:sim-Gi}
\end{figure}

We now consider the NFC--Gilad model \eqref{eqn:Gilad} with parameter values \eqref{pars:Gilad}. {  We recall from Section \ref{ss:Gi} that
\begin{equation*}
    P_1 = 0.381>0,\qquad P_2 = -0.207<0, \qquad P_3 = -0.575<0, \qquad P_4 = 0.818>0
\end{equation*}
and so we expect localised anti-phase gaps to emerge as the precipitation increases from the Turing point $\mu_*=1.635$. Hence, we simulate the time evolution of our initial guess \eqref{InitialGuess}---for hexagons, squares and pentagons---in the parameter region $\mu>\mu_*$.}

The Turing bifurcation has an associated wave number $k\approx 0.333$, and so we compute our simulations in a square box of length $378$. The simulations for our three examples are presented in Figure~\ref{fig:sim-Gi} for the vegetation density $u$; solutions are again found to be anti-phase, and so we {  mostly} do not plot the water density $v$ since its spatial structure is identical to $u$. 

Much like with the logistic Klausmeier model \eqref{eqn:extK}, individual gaps in the NFC--Gilad model grow in width as time evolves. However, unlike the logistic Klausmeier model, these gaps do not coalesce with their neighbours. As a result, we observe patterns with non-uniform elliptic gaps that retain a mosaic-like structure. We note that the interface between the pattern and the surrounding vegetated state is much more circular for \eqref{eqn:Gilad} than in our previous models, even for the square and pentagon examples.  

We note that all three examples exhibit a pseudo-hexagonal lattice structure away from their core, which might hint towards the dominance of hexagonal patterns in field observations and experiments. One could then think of the resultant square and pentagon patterns as versions of the hexagonal example with defects at their core. However, it is worth noting that the hexagonal packing observed in the square and pentagon examples is not symmetry breaking; the dihedral symmetry of each example ($\mathbb{D}_{4}$ for squares, $\mathbb{D}_{5}$ for pentagons) is preserved throughout, except for when the pentagon example becomes large enough to experience boundary effects.

\,
\subsubsection{Von Hardenberg - Turing point 1}
\begin{figure}[t!]
    \centering
    \includegraphics[width=0.9\linewidth]{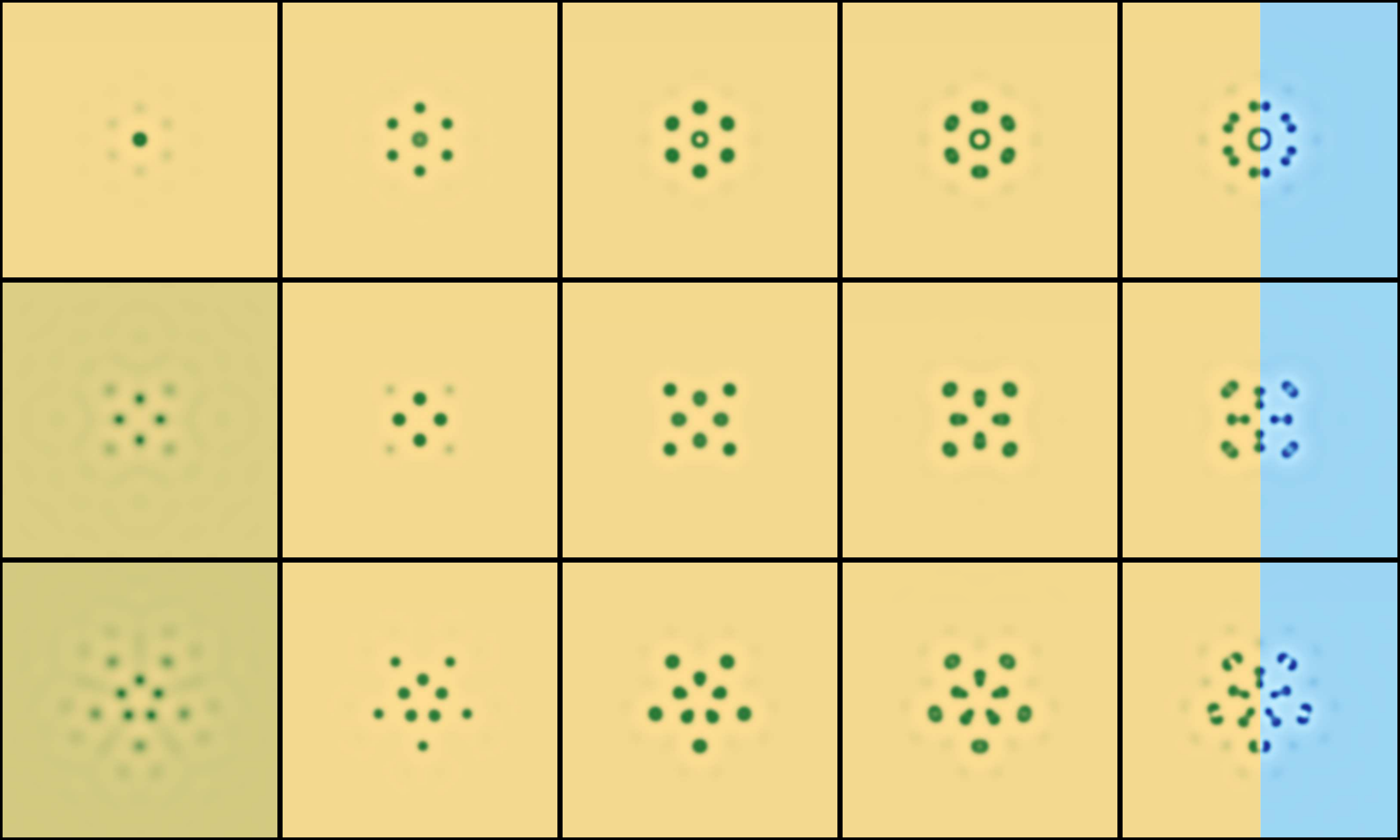}
    \caption{Time evolution of vegetation density $u$ in the von Hardenberg model \eqref{eqn:vonH} with parameter values \eqref{pars:vonH} near its first Turing point. The dark green represents a vegetated state ({  $u\approx0.35$}), while the lighter yellow represents a bare state ({  $u\approx0$}). The rows are the same as Figure~\ref{fig:sim-KG}, the columns correspond to the points in time $t=200, 300, 400, 500$ and $600$. {  The final column is divided into vegetation density (left) and water density (right), where darker and lighter blue represents higher and lower water density, respectively. Videos of each simulation can be found at \cite{Hill2023Github-Veg}.}}
    \label{fig:sim-vH-1}
\end{figure}

We consider the first Turing bifurcation \eqref{Turing:vH-1} in the von Hardenberg model \eqref{eqn:vonH} with parameter values \eqref{pars:vonH}. {  We recall from Section \ref{ss:vH} that, for this bifurcation point,
\begin{equation*}
    P_1 = -0.384<0,\qquad P_2 = 1.707>0, \qquad P_3 = 0.843>0, \qquad P_4 = 0.001>0
\end{equation*}
and so we expect localised in-phase peaks to emerge as the precipitation decreases from the Turing point $\mu_*=0.169$. Hence, we simulate the time evolution of our initial guess \eqref{InitialGuess}---for hexagons, squares and pentagons---in the parameter region $\mu<\mu_*$.}

This bifurcation has an associated wave number $k\approx 0.106$, and so we compute our simulations in a square box of length $1190$. However, the width of the patterns is considerably smaller than this domain, and so we present our results in the smaller square box of length $595$. The water density $v$ is found to be in-phase, as predicted in Section~\ref{s:Theory}, and shares an identical spatial structure to the vegetation density $u$. Hence, we present our simulations of $u$ in Figure~\ref{fig:sim-vH-1}. The evolution of these patterns occurs at a slower rate than in previous models, and so we present our simulations at later times, as detailed in Figure~\ref{fig:sim-vH-1}. 

The initial guess \eqref{InitialGuess} is very small for this Turing bifurcation, and so it can be difficult to capture numerically. As such, in order to observe localised dihedral patterns in this model, we take a larger choice for the constant $C$ in \eqref{InitialGuess}. Then, for each example solutions converge to localised peaks spread out in dihedral arrangements. The emergence of new peaks is much slower than in the previous models; individual spots instead grow in width until their core destabilises, causing a transition into a ring or two smaller spots. This behaviour is reminiscent of localised spikes in singularly perturbed reaction--diffusion systems, such as  \cite{Kolokolnikov2022ring}, rather than the invading patterns seen in our other models. We note that we observe these spot-to-ring transitions occur even though localised ring-type patterns do not emerge from the uniform state.

Other than the central spot in the hexagon example, every spot is observed transitioning into two smaller spots. The way that each spot splits appears to be very structured; the outer ring of spots appears to split in the azimuthal (i.e. angular) direction, whereas the inner ring of spots appears to split in the radial direction. Since the outer layer of spots also appears to split symmetrically, we note that the overall dihedral symmetry of the system seems to be preserved by these splittings.

\subsubsection{Von Hardenberg - Turing point 2}
\begin{figure}[t!]
    \centering
    \includegraphics[width=0.9\linewidth]{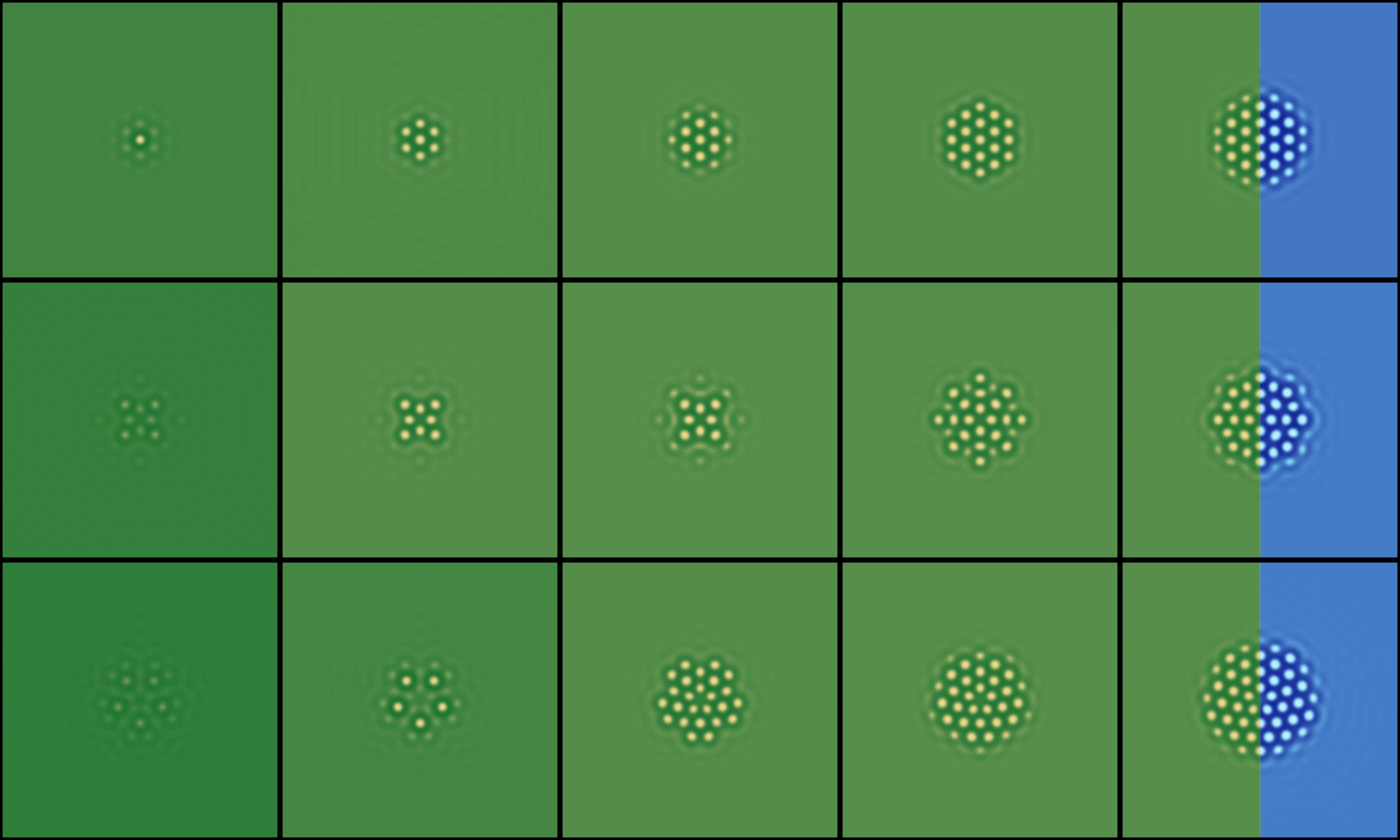}
    \caption{Time evolution of vegetation density $u$ in the von Hardenberg model \eqref{eqn:vonH} with parameter values \eqref{pars:vonH} near its second Turing point. The dark green represents a vegetated state ({  $u\approx0.35$}), while the lighter yellow represents a bare state ({  $u\approx0$}). The rows and columns are the same as Figure~\ref{fig:sim-vH-1}. {  The final column is divided into vegetation density (left) and water density (right), where darker and lighter blue represents higher and lower water density, respectively. Videos of each simulation can be found at \cite{Hill2023Github-Veg}.}}
    \label{fig:sim-vH-2}
\end{figure}

Our final simulations are for the second Turing point \eqref{Turing:vH-2} in the von Hardenberg model \eqref{eqn:vonH} with parameter values \eqref{pars:vonH}. {  We recall from Section \ref{ss:vH} that, for this bifurcation,
\begin{equation*}
    P_1 = 0.217>0, \qquad P_2 = 2.578>0, \qquad P_3 = -1.512<0, \qquad P_4 = 0.014>0
\end{equation*}
and so we expect localised in-phase gaps to emerge as the precipitation increases from the Turing point $\mu_*=0.414$. Hence, we simulate the time evolution of our initial guess \eqref{InitialGuess}---for hexagons, squares and pentagons---in the parameter region $\mu>\mu_*$.}

This bifurcation has an associated wave number $k\approx 0.201$, and so we compute our simulations in a square box of length $612$. Again, we find the water density $v$ to be in-phase, as predicted in Section~\ref{s:Theory}, with an identical spatial structure to the vegetation density $u$; we present our simulations of $u$ in Figure~\ref{fig:sim-vH-2}. Again, the time evolution of our patterns is slower than the other models, and so we present our simulations for slightly later times. 

In each example solutions are comprised of circular gaps with a uniform size. The hexagon example results in a perfect hexagonal lattice, exhibiting the same behaviour as in the Klausmeier--Gray--Scott model and prototypical pattern-forming systems like the Swift--Hohenberg equation \cite{lloyd2008localized}. The square and pentagon examples grow similarly to the Klausmeier--Gray--Scott \eqref{eqn:KGS}, but new gaps exhibit hexagonal packing like in the NFC--Gilad models \eqref{eqn:Gilad}. Then, each example evolves into a perfect hexagonal lattice with possible defects; these defects then form grain boundaries throughout the localised pattern, such as the horizontal and vertical lines in the square example. As with the previous models, even though the solutions are all evolving into distorted hexagonal lattices, the dihedral symmetry of the solution is preserved.
\,
\section{Discussion}
In this work, we have provided a step-by-step guide on how to convert a two-component reaction--diffusion system into its local representation near a Turing instability for a class of vegetation models. Then one can find localised dihedral patterns bifurcating from the Turing instability and obtain four quantities that encode the qualitative behaviour of these localised solutions and how they bifurcate. We found that localised spot A-type patterns emerge generically from a Turing instability, not requiring any subcriticality conditions like localised one dimensional and ring-type patterns.

We then considered four established models for dryland vegetation and predicted the emergence of localised dihedral patterns. For the simple Klausmeier--Gray--Scott model \eqref{eqn:KGS}, we were able to derive explicit forms for each qualitative predictor in terms of the parameters of the model, and hence better understand the qualitative robustness of localised dihedral solutions with regards to our choice of parameters. We found that each model exhibits spot A-type patterns bifurcating from a Turing point, but ring-type patterns do not emerge for the given choices of parameters from the literature. However, in the Klausmeier--Gray--Scott model we saw that there exists a small parameter region in which ring-type patterns emerge from the uniform state, which is also observed for the von Hardenberg model in \cite[Figure~15]{Dawes2016}. We then supported our predictions through numerical simulations of each vegetation model, where we observed that the generic spot A-type patterns evolve into dramatically different structures depending on the choice of model. 

We note that the dihedral symmetry of each solution is very robust in our simulations, even though the Cartesian discretisation of our spatial domain should actively disfavour the rotational symmetry of our patterns. One possible reason for this might be that the emergence of each peak or gap is strongly influenced by its neighbours. Then, if we begin with a collection of peaks or gaps arranged with dihedral symmetry, the neighbourly effects on a new peak or gap also possess dihedral symmetry, thus enforcing the dihedral structure on the new peaks and gaps. Another possible reason might be that our initial guess imparts a dihedral symmetry on the whole domain which the localisation just makes extremely small away from the origin. Then, as the uniform state destabilises, the already present dihedral oscillations grow in amplitude, resulting in a pattern with dihedral symmetry.

Although our solutions maintain their dihedral symmetry as time evolves, we observe significant differences in the resultant patterns for different vegetation models. The Klausmeier--Gray--Scott \eqref{eqn:KGS}, NFC--Gilad \eqref{eqn:Gilad} and von Hardenberg models \eqref{eqn:vonH} all exhibit similar growth behaviour in Figures~\ref{fig:sim-KG}, \ref{fig:sim-Gi} \& \ref{fig:sim-vH-2}, but with significantly different final patterns. For the Klausmeier--Gray--Scott model in Figure~\ref{fig:sim-KG} there is no clear packing structure as each pattern grows larger, whereas solutions the NFC--Gilad model in Figure~\ref{fig:sim-Gi} exhibit nonuniform hexagonal packing and solutions for the von Hardenberg model in Figure~\ref{fig:sim-vH-2} exhibit uniform hexagonal packing. We note that we cannot reasonably compare the speed of growth of the patterns in each model, as each model considered in this work is nondimensionalised with a rescaled time coordinate. For example, for the fixed parameter values given, $t=500$ in our simulations corresponds to approximately 50 years in the NFC--Gilad model but 125 years in the logistic Klausmeier model.

 We have provided novel tools for the design of phenomenological reaction--diffusion models; by computing the quantities $P_1,P_2,P_3,P_4$, one is able to quickly check the qualitative behaviour of solutions of a prospective model and compare them against the observed phenomena. {  While these predictors are derived for localised patterns, they also have applications in general pattern formation. Like spot A-type patterns, domain-covering hexagons also possess a fixed phase and polarity---as predicted by $P_2$ and $P_3$, respectively---while $P_4$ determines whether unstable stripes bifurcate in the direction determined by $P_1$. In order to model the emergence of a particular type of pattern, one would then require a particular type of Turing point determined by the signs of each $P_i$; this is summarised in Table~\ref{Table:1}. 
\begin{table}[ht!]
     \centering
    \begin{tabular}{|c|c|c|}
         \hline Patterns & Example & Requirement \\
         \hline
         In-phase spots & \cite{Yizhaq2014heterogeneous} & $P_2>0$,\quad  $P_3>0$\\
         Anti-phase spots & \cite{Yizhaq2014heterogeneous} & $P_2<0$,\quad  $P_3>0$\\
         In-phase gaps & \cite{Getzin2016FairyCircles} & $P_2>0$,\quad  $P_3<0$\\
         Anti-phase gaps & \cite{Getzin2016FairyCircles} & $P_2<0$,\quad  $P_3<0$\\
         Rings and (unstable) stripes & \cite{sheffer2007plants} & $P_4<0$\\
         \hline
    \end{tabular}
    \caption{A summary of the different qualitative patterns emerging from a Turing point and their required values for $P_2, P_3$ and $P_4$.}
    \label{Table:1}
\end{table}
}

For example, if you were looking to model localised in-phase gaps (such as fairy circles observed in Australia \cite{Getzin2016FairyCircles}), then any prospective reaction--diffusion model must have a Turing instability with $P_2>0$, $P_3<0$. {  Such a model would not be appropriate for the fairy circles observed in Namibia, however, which manifest as anti-phase gaps and thus require a model with $P_2<0$, $P_3<0$. Furthermore, we note that any vegetation model that exhibits a transition from spots to gaps (or vice versa), such as the von Hardenberg model~\cite{Gowda2014}, must thus possess two distinct Turing points with different signs for $P_3$.} While these properties alone are not sufficient to successfully model complex phenomena such as vegetation patterns, they can provide helpful conditions to be used in conjunction with ecological modelling. We also note that the spot A-type patterns emerge generically from a Turing bifurcation, and so the mere existence of localised patterns is not necessarily a good indicator of whether or not a model is appropriate, since if a model can exhibit two-dimensional periodic patterns it can also exhibit localised dihedral patterns.

 Beyond the design of new models, our work also provides some insight into the formation of vegetation patterns. In our simulations, localised patterns represent a transition between two states; a stable patterned state and an unstable uniform state. The patterned state emerges at the origin and begins to invade the uniform state in a radial direction, as seen prototypical pattern-forming systems \cite{Lloyd2019InvasionFronts,Lloyd2021HexInvasionFronts}; this is particularly serious in the logistic Klausmeier model \eqref{eqn:extK} where the invading patterned state collapses to the bare state, resulting in the desertification of the local environment. Since localised dihedral patterns lead to the further destabilisation of the uniform vegetated state, they may serve the same role as periodic patterns in predicting or mitigating tipping points and desertification effects in semi-arid environments \cite{Rietkerk2004catastrophe,Rietkerk2021Tipping}.

We briefly comment on the popular phenomenon of fairy circles and their relation to the localised gaps discussed in this work. In \cite{Getzin2021Fairy}, Getzin et al. define fairy circles as satisfying the following conditions; (a) they must be gaps in a vegetated state, that (b) can form periodic structures, and (c) are not rings surrounded by a bare state. We have found that, for any Turing bifurcation with $P_3<0$, there are localised dihedral patterns that satisfy these criteria, as seen in Figures~\ref{fig:sim-KG}, \ref{fig:sim-Kl}, \ref{fig:sim-Gi} \& \ref{fig:sim-vH-2}. Furthermore, we have found that these gap solutions are distinct from ring-type patterns found in Theorem~\ref{thm:Ring}, which emerge in much more restrictive parameter regimes.

{  We note that there are numerous possible directions in which to numerically explore localised vegetation patterns, many of which we do not consider here. The theory in Section~\ref{s:Theory} describes the bifurcation of localised patterns from a Turing instability and so a natural first step would be to use numerical continuation schemes to characterise the bifurcation curves of these patterns in more detail. One could adapt the pseudo-spectral codes of \cite{Hill2023DihedralPatch} to each vegetation model, using a truncated angular Fourier expansion and finite difference to solve and continue solutions with a given $\mathbb{D}_{m}$ dihedral symmetry. However, determining the nonlinearities of each model projected onto angular Fourier modes remains a difficult task, and each bifurcation curve would be restricted to the dihedral symmetry group $\mathbb{D}_{m}$, thus preventing symmetry-breaking instabilities. One could instead use a PDE solver such as \verb!pde2path! \cite{uecker_wetzel_rademacher_2014}, but the large domains required for vegetation models would make any continuation study computationally intensive. It would be interesting to consider the bifurcation structure of these localised patterns in a future numerical study. Furthermore, it would also be interesting to explore the sensitivity of our localised vegetation patterns to the domain boundary. The size of our spatial domain and choice of boundary conditions can both cause nonlinear effects in the time evolution of vegetation patterns; for localised patterns in sufficiently large domains these effects are negligible, but become more important as the vegetation patches grow in width. We again leave this topic for future study.}

A possible extension of this work would be to consider nonlocal equations or reaction--diffusion systems with more than 2 components; there are several examples of vegetation models of these forms, including in \cite{Gilad2004Engineers} and \cite{Rietkerk2002}. However, such problems would require some kind of model reduction, either formally or rigorously, before one could apply the techniques presented in this work. It would also be interesting to try and investigate more complicated two-dimensional localised patterns, such as localised labyrinthine patterns observed numerically \cite{Clerc2021Labyrinthine}. Of course, the theory presented in this work goes beyond the study of vegetation patterns, and could be applied to any two-component reaction--diffusion system that undergoes a Turing bifurcation. Some interesting examples would include the Lugiato--Lefever equation in nonlinear optics, and the Gross-Pitaevskii (GP) equation for Bose-Einstein solitons in quantum mechanics. 

This work also highlights the need for more topological measures in order to compare localised planar patterns. We discussed the qualitative differences between our numerical simulations, but we are not currently equipped with the tools to easily compare them quantitatively. There have been recent developments in this direction for studying nearly hexagonal lattices \cite{Subramanian2021Grain,Motta2018Measures} which could provide further insight into the dynamics observed in our simulations.

\begin{Data}
    All data used to produce the results of this paper is openly available at the following URL: \href{https://github.com/Dan-Hill95/Dihedral_Vegetation}{https://github.com/Dan-Hill95/Dihedral{\textunderscore}Vegetation} \cite{Hill2023Github-Veg}.
\end{Data}

\begin{Acknowledgment}  
   {  The author gratefully acknowledges support from the Alexander von Humboldt Foundation. They} would {  also} like to thank Josephine Solowiej--Wedderburn, Nicolas Verschueren, Nils Gutheil, Rami Ahmad, David Lloyd and Tom Bridges for their feedback on an earlier draft of this manuscript, {  as well as the anonymous referees for their helpful comments and suggestions.}
\end{Acknowledgment}

\appendix

\setcounter{equation}{0}
\renewcommand\theequation{\Alph{section}.\arabic{equation}}

\section{Analysis of the Klausmeier--Gray--Scott model}\label{app:KGS-analysis}
We now explicitly analyse the Klausmeier--Gray--Scott model \eqref{eqn:KGS}
\begin{equation*}
    \begin{split}
        u_t &= \Delta u + v u^2 - m u,\\
        v_t &= \delta_v \Delta v + \mu - v - v u^2.
    \end{split}
\end{equation*}
We begin by finding uniform steady states, which satisfy 
\begin{equation*}
    \begin{split}
        v = \frac{\mu}{1+u^2}, \qquad \qquad 0 = \left[u^2 - \left(\frac{\mu}{m}\right) u + 1\right]m u.
    \end{split}
\end{equation*}
In particular, we obtain a trivial solution $u_0 = 0$, $v_0 = \mu$, as well as two non-trivial solutions
\begin{equation*}
    \begin{split}
        u_\pm(\mu) = \left(\frac{\mu}{2m}\right) \pm \sqrt{\left(\frac{\mu}{2m}\right)^2 - 1} \qquad  v_{\pm}(\mu) = \frac{m}{u_{\pm}(\mu)},
    \end{split}
\end{equation*}
which exist for $\mu \geq 2m$. The trivial steady state never undergoes a Turing instability as $\mu$ varies, and so we linearise \eqref{eqn:KGS} about $(u_\pm, v_\pm)$; that is, we define $u(t,\mathbf{x}) = u_\pm(\mu) + U(t,\mathbf{x})$, $v(t,\mathbf{x}) = \frac{m}{u_\pm(\mu)} + V(t,\mathbf{x})$, such that \eqref{eqn:KGS} becomes
\begin{equation}\label{eqn:KGS-pert}
    \begin{split}
        0 &= \Delta U + m U  + u_\pm(\mu)^2 V + \frac{m}{u_\pm(\mu)} U^2 + 2u_\pm(\mu) U V + U^2 V,\\
        0 &= \Delta V - \frac{1}{\delta_v}\left[2m U  + (1 + u_\pm(\mu)^2) V + \frac{m}{u_\pm(\mu)}U^2 + 2 u_\pm(\mu) UV + U^2V\right].
    \end{split}
\end{equation}
We now want to identify any values of $\mu$ such that \eqref{eqn:KGS-pert} undergoes a Turing instability. First, we write \eqref{eqn:KGS-pert} in the vector form of \eqref{eqn:Taylor-mu}, so that
\begin{equation*}
    \mathbf{0} = \Delta \mathbf{U} - \mathbf{M}(\mu)\mathbf{U} - \mathbf{F}_{2}(\mathbf{U};\mu) - \mathbf{F}_{3}(\mathbf{U};\mu),
\end{equation*}
with 
\begin{equation*}
    \mathbf{M}(\mu) = \begin{pmatrix}
        - m & - u_\pm^2(\mu) \\
        \tfrac{2m}{\delta_v } & \tfrac{1}{\delta_v }\left(1 + u_\pm^2(\mu)\right)
    \end{pmatrix},
\end{equation*}
and
\begin{equation*}
    \mathbf{F}_{2}(\mathbf{U};\mu) = \begin{pmatrix}
        -\frac{m}{u_\pm(\mu)} U^2 - 2u_\pm(\mu) U V\\
        \frac{1}{\delta_v}\left[\frac{m}{u_\pm(\mu)} U^2 + 2u_\pm(\mu) U V\right]
    \end{pmatrix}, \qquad \mathbf{F}_{3}(\mathbf{U};\mu) = \begin{pmatrix}
        - U^2 V\\
        \frac{1}{\delta_v} U^2 V
    \end{pmatrix}.
\end{equation*}
A general two-dimensional matrix $\mathbf{A}$ has a double eigenvalue of $\lambda = \frac{1}{2}\mathrm{tr}\,\mathbf{A}$ if and only if 
\begin{equation*}\begin{split}
    (\mathrm{tr}\,\mathbf{A})^2 - 4\,\det\mathbf{A} ={}& 0,
\end{split}
\end{equation*}
where $\mathrm{tr}\,\mathbf{A}$ and $\det\,\mathbf{A}$ denote the trace and determinant of the matrix $\mathbf{A}$. For the above matrix $\mathbf{M}(\mu)$, we obtain 
\begin{equation*}\begin{split}
    (\mathrm{tr}\,\mathbf{M})^2 - 4\,\det\mathbf{M} ={}& m^2 - 2m\tfrac{1}{\delta_v }\left(1 + u_\pm^2\right) + \tfrac{1}{\delta^2_v }\left(1 + u_\pm^2\right)^2 - 4\tfrac{m}{\delta_v }\left(u_\pm^2-1\right),\\
    ={}&\tfrac{1}{\delta_v^2}\left[ (\delta_v m + 1)^2 + 2 (1 - 3(\delta_v m))u_\pm^2 + u_\pm^4\right],
\end{split}
\end{equation*}
and hence any repeated eigenvalues must satisfy
\begin{equation*}
    u_*^{\pm} = \sqrt{3(\delta_v m) - 1 \pm \sqrt{(3(\delta_v m)-1)^2 - (\delta_v m + 1)^2}} = \sqrt{3(\delta_v m) - 1 \pm \sqrt{8(\delta_v m)^2 - 8(\delta_v m)}}
\end{equation*}
which implies that
\begin{equation*}
    \frac{\mu_{*}^\pm}{m} = \frac{1+(u_*^\pm)^2}{u_*^\pm} = \frac{3(\delta_v m) \pm \sqrt{8(\delta_v m)^2 - 8(\delta_v m)}}{\sqrt{3(\delta_v m) - 1 \pm \sqrt{8(\delta_v m)^2 - 8(\delta_v m)}}}.
\end{equation*}
Computing the eigenvalue $\lambda = \frac{1}{2}\mathrm{tr}\,\mathbf{M}(\mu_*^{\pm})$, we find that
\begin{equation*}
    \lambda^{\pm} = \tfrac{1}{2\delta_v}\left[ (1+(u_*^\pm)^2) - (\delta_v m)\right] = \tfrac{1}{\delta_v}\left[ (\delta_v m) \pm \sqrt{(\delta_v m)^2 + [(\delta_v m)- 2](\delta_v m)}\right].
\end{equation*}

Hence, we have found two bifurcations with repeated roots; however, $\lambda^{+}$ (corresponding to $u_*^+$) is positive for all $(\delta_v m)>0$ and $\lambda^{-}$ (corresponding to $u_*^-$) is only negative if $(\delta_v m)>2$.  

Let us briefly make a couple of observations. Firstly, we note that both points $u_{*}^{\pm}$ lie on the upper branch of the nontrivial uniform steady state, such that $u_*^\pm = u_{+}(\mu_*^\pm)\geq1$. Secondly, we note that the upper point $u_*^+$ is called a Belyakov--Devaney point, since it always possesses a positive repeated eigenvalue. These points are related to spot to spike transitions; see \cite{alSaadi2021GrayScott} for an example of such transitions.  

We continue by assuming that $(\delta_v m)>2$ so that the equilibrium
\begin{equation*}
    (u,v) = \left( u_{+}(\mu), \frac{m}{u_{+}(\mu)}\right),
\end{equation*}
undergoes a Turing bifurcation at $\mu=\mu_*^-$, $u=u_*:=\sqrt{3(\delta_v m) - 1 - \sqrt{8(\delta_v m)^2 - 8(\delta_v m)}}$ with wave number $k>0$, defined by
\begin{equation*}
    k = \sqrt{\tfrac{1}{2\delta_v}\left[(\delta_v m) - (1 + u_*^2)\right]} = \sqrt{\tfrac{1}{\delta_v}\left(\sqrt{2(\delta_v m)^2 - 2(\delta_v m)} - (\delta_v m)\right)}.
\end{equation*}
Taking $\mu = \mu_*^- + \varepsilon$, with $|\varepsilon|\ll1$, we note that 
\begin{equation*}
    \begin{split}
        u_+(\mu_*^- + \varepsilon) = u_* + \frac{1}{m}\left(\frac{u_*^2}{u_*^2-1}\right)\varepsilon + \mathcal{O}(\varepsilon^2).
    \end{split}
\end{equation*}
Then, we can now express \eqref{eqn:KGS-pert} in the form of \eqref{eqn:Taylor-eps}, 
\begin{equation*}
    \begin{split}
        0 &= \Delta \mathbf{U} - \mathbf{M}_{1}\mathbf{U} -\varepsilon\mathbf{M}_{2}\mathbf{U} - \mathbf{Q}(\mathbf{U},\mathbf{U}) - \mathbf{C}(\mathbf{U},\mathbf{U},\mathbf{U}) + \mathcal{O}\left(|\varepsilon|\,|\mathbf{U}|^2 + |\varepsilon|^2\,|\mathbf{U}|\right),\\
    \end{split}
\end{equation*}
where 
\begin{equation*}\begin{split}
    \mathbf{M}_{1} :={}& \begin{pmatrix} -m & -u_*^2 \\
    \frac{2m}{\delta_v} & \frac{1}{\delta_v}(1+u_*^2)
    \end{pmatrix} = \begin{pmatrix} -m & -u_*^2 \\
    \frac{(m-k^2)^2}{u_*^2} & m - 2 k^2
    \end{pmatrix}, \qquad  \mathbf{M}_{2} :={} \begin{pmatrix} 0 & -\left(\frac{2 u_*^3}{m(u_*^2-1)}\right) \\
    0 & \left(\frac{u_*(m-k^2)^2}{m^2(u_*^2-1)}\right)
    \end{pmatrix},\\
\end{split}\end{equation*}
and
\begin{equation*}\begin{split}
    \mathbf{Q}(\mathbf{U}_1,\mathbf{U}_2) :={}& \left[\frac{m}{u_*} U_1U_2 + u_* (U_2V_1 + U_1 V_2)\right]\begin{pmatrix}
        -1\\
        \frac{(m-k^2)^2}{2m u_*^2}
    \end{pmatrix}, \\
        \mathbf{C}(\mathbf{U}_1,\mathbf{U}_2,\mathbf{U}_3) :={}& \tfrac{1}{3}\left[U_1U_2V_3 + U_1U_3V_2 + U_2U_3V_1\right]\begin{pmatrix}
        -1\\
        \frac{(m-k^2)^2}{2m u_*^2}
    \end{pmatrix}.\\
\end{split}\end{equation*}
Here, we have used the relation
\begin{equation*}
    u_*^2 = \delta_v \frac{(m-k^2)^2}{2m}
\end{equation*}
in order to write down $\mathbf{M}_1,\mathbf{M}_2,\mathbf{Q}$ and $\mathbf{C}$ in terms of $m,u_*$ and $k$. We introduce the generalised eigenvectors $\hat{U}_{0},\hat{U}_{1}$ and their dual basis $\hat{U}_{0}^{*},\hat{U}_{1}^{*}$, 
\begin{equation*}
    \hat{U}_0 = \begin{pmatrix}
        -u_*^2 \\ m  - k^2
    \end{pmatrix}, \qquad \hat{U}_1 = \begin{pmatrix}
        0\\ k^2
    \end{pmatrix}, \qquad \hat{U}^*_0 = \begin{pmatrix}
        -\frac{1}{u_*^2} \\ 0
    \end{pmatrix},  \qquad\hat{U}^*_1 = \begin{pmatrix}
        \frac{m - k^2}{k^2 u_*^2} \\ \frac{1}{k^2}
    \end{pmatrix},
\end{equation*}
so that
\begin{equation*}\begin{split}
    (\mathbf{M}_1 + k^2\mathbbm{1})\hat{U}_0 ={}& \begin{pmatrix} -(m-k^2) & -u_*^2 \\
    \frac{(m-k^2)^2}{u_*^2} & m - k^2
    \end{pmatrix}\begin{pmatrix}
        -u_*^2 \\ m  - k^2
    \end{pmatrix}= \begin{pmatrix}
        (m-k^2)u_*^2 - u_*^2(m-k^2)\\
        -\frac{(m-k^2)^2}{u_*^2}u_*^2 + (m-k^2)^2
    \end{pmatrix}=\mathbf{0},\\
    (\mathbf{M}_1 + k^2\mathbbm{1})\hat{U}_1 ={}& \begin{pmatrix} -(m-k^2) & -u_*^2 \\
    \frac{(m-k^2)^2}{u_*^2} & m - k^2
    \end{pmatrix}\begin{pmatrix}
        0\\ k^2
    \end{pmatrix} = \begin{pmatrix}
        -k^2 u_*^2 \\ k^2(m-k^2)
    \end{pmatrix} = k^2 \hat{U}_0
\end{split}
\end{equation*}
and $\hat{U}_{i}\cdot\hat{U}_{j} = \delta_{i,j}$. Finally, we can determine the explicit forms of the qualitative predictors $P_j$; starting with $P_1$, we see that
\begin{equation*}
    \begin{split}
    P_1 =c_0 ={}& \hat{U}_1^* \cdot \left(-\tfrac{1}{4}\mathbf{M}_2\hat{U}_0\right) \\
 ={}& -\frac{1}{4} \left(\frac{2 u_*^3}{m(u_*^2-1)}\right)\begin{pmatrix}
        \frac{m - k^2}{k^2 u_*^2} \\ \frac{1}{k^2}
    \end{pmatrix}\cdot\begin{pmatrix} 0 & -1 \\
    0 & \frac{(m-k^2)^2}{2m u_*^2}
    \end{pmatrix}\begin{pmatrix}
        -u_*^2 \\ m  - k^2
    \end{pmatrix}\\
    ={}& -\frac{1}{4} \left(\frac{2 u_*^3(m  - k^2)}{m(u_*^2-1)}\right)\begin{pmatrix}
        \frac{m - k^2}{k^2 u_*^2} \\ \frac{1}{k^2}
    \end{pmatrix}\cdot\begin{pmatrix} -1 \\
    \frac{(m-k^2)^2}{2m u_*^2}
    \end{pmatrix}\\
    ={}& \left(\frac{u_*(m  - k^2)^2(m+k^2)}{4 k^2 m^2 (u_*^2-1)}\right).\\
    \end{split}
\end{equation*}
Since $u_* = u_{+}(\mu_*^-)$, we know that $u_*^2 - 1 > 0$, and so $P_1>0$ for all parameter values. Turning our attention to $P_2$, we find that
\begin{equation*}
    P_2 = \frac{[U_0]_2}{[U_0]_1} = -\left(\frac{m-k^2}{u_*^2}\right), \\
\end{equation*}
and so $P_2<0$ for all parameter values, since $m-k^2 = \frac{1}{2}[ m + \frac{1}{\delta_v}(1 + u_*^2)] > 0$. Before computing $P_3$ and $P_4$, we compute the following nonlinear terms
\begin{equation*}
\begin{split}
    Q_{00} ={}& \mathbf{Q}(\hat{U}_{0},\hat{U}_{0})= -u_*(m-2k^2)\begin{pmatrix}
        -u_*^2 \\ \frac{(m-k^2)^2}{2m}
    \end{pmatrix},\\
    Q_{01} ={}& \mathbf{Q}(\hat{U}_{0},\hat{U}_{1}) = - k^2 u_*\begin{pmatrix}
        -u_*^2 \\ \frac{(m-k^2)^2}{2m}
    \end{pmatrix},\\
    C_{000} ={}& \mathbf{C}(\hat{U}_{0},\hat{U}_{0},\hat{U}_{0}) = u_*^2 (m-k^2)\begin{pmatrix}
        -u_*^2 \\ \frac{(m-k^2)^2}{2m}
    \end{pmatrix}.\\
\end{split}
\end{equation*}
Then, we see that
\begin{equation*}
\begin{split}
    P_3 ={}& \frac{[U_0]_{1}}{\hat{U}^*_1\cdot Q_{00}} =-\left(\frac{2 u_* m k^2}{(m-2k^2)(m - k^2)(m+k^2)}\right),
\end{split}
\end{equation*}
and so $P_3<0$ for all parameter values, since $m - 2k^2 = \frac{1}{\delta_v}(1+u_*^2)>0$. Finally, we obtain
\begin{equation*}
\begin{split}
   P_4 =  c_3 ={}& -\frac{5}{6}\left(\hat{U}_0^*\cdot Q_{00}\right)\left(\hat{U}_1^*\cdot Q_{00}\right) -\frac{5}{6}\left(\hat{U}_1^*\cdot Q_{01}\right)\left(\hat{U}_1^*\cdot Q_{00}\right) -\frac{19}{18}\left(\hat{U}_1^*\cdot Q_{00}\right)^2 - \frac{3}{4} \left(\hat{U}_1^*\cdot C_{000}\right) ,\\
      ={}& \frac{5}{6}\left(\frac{u_*^2(m-2k^2)^2(m - k^2)(m+k^2)}{2 m k^2}\right)  -\frac{5}{6}\left(\frac{u_*^2(m-2k^2)(m - k^2)^2(m+k^2)^2}{4 m^2 k^2}\right)\\
    {}&\qquad -\frac{19}{18}\left(\frac{u_*^2(m-2k^2)^2(m - k^2)^2(m+k^2)^2}{4 m^2 k^4}\right) + \frac{3}{4} \left(\frac{u_*^2(m-k^2)^2(m+k^2)}{2mk^2}\right).\\
\end{split}
\end{equation*}
It is not clear for each $m>0$, $\delta_v>\frac{2}{m}$ whether $P_4>0$ or $P_4<0$, and so we numerically plot the sign of $P_4$ in Figure~\ref{fig:K-GS}(a).

\bibliographystyle{abbrv}
\bibliography{Reference}
\end{document}